\documentclass[11pt, reqno]{amsart}

\usepackage{amsmath}
\usepackage{fullpage}
\usepackage{amssymb}
\usepackage{amsthm}
\usepackage{bbm}
\usepackage{subfigure}
\usepackage{float}
\usepackage{enumitem}
\usepackage{graphicx}
\usepackage{color}
\usepackage{epstopdf}
\usepackage{enumerate}
\graphicspath{{../fig/}}

\newtheorem{lemma}{Lemma}[section]
\newtheorem{theorem}[lemma]{Theorem}
\newtheorem{proposition}[lemma]{Proposition}

\newtheorem{corollary}[lemma]{Corollary}
\newtheorem{definition}[lemma]{Definition}
\allowdisplaybreaks

\newcommand{\fe}{{\rm e}}

\DeclareMathOperator{\im}{Im}
\numberwithin{equation}{section}
\newcommand{\msc}[2][2000]{
  \let\@oldtitle\@title
  \gdef\@title{\@oldtitle\footnotetext{#1 \emph{Mathematics subject
        classification.} #2}}
}

\title[A low regularity EI for the dNLS]{A low regularity exponential-type integrator for the derivative nonlinear Schr\"odinger equation}

\author{Lun Ji}
\address{Department of Applied Mathematics, the Hong Kong Polytechnic University, Hung Hom, Hong Kong (L.~Ji)}
\email{lun-422.ji@polyu.edu.hk}

\author{Hang Li}
\address{LJLL (UMR 7598), Sorbonne Universit\'e, UPMC, 4 place Jussieu, 75005 Paris, France (H.~Li)}
\email{li727263992@gmail.com}

\author{Alexander Ostermann}
\address{Digital Science Center, Universit\"at Innsbruck, 6020 Innsbruck, Austria (A.~Ostermann)}
\email{alexander.ostermann@uibk.ac.at}


\thanks {L. Ji is partially supported by the Research Grants Council of Hong Kong (grant No.~15306123). H. Li is supported by the European Research Council (ERC) under the European Union's Horizon 2020 research and innovation program (grant agreement No.~850941), as well as by a postdoctoral fellowship from the Fondation Sciences Math\'ematiques de Paris (FSMP)}

\begin{document}

\begin{abstract}
In this work, we present a first-order unfiltered exponential integrator for the one-dimensional derivative nonlinear Schr\"odinger equation with low regularity. Our analysis shows that for any $s>\frac12$, the method converges with first-order in $H^s(\mathbb{T})$ for initial data $u_0\in H^{s+1}(\mathbb{T})$. Moreover, we constructed a symmetrized version of this method that performs better in terms of both global error and conservation behavior. To the best of our knowledge, these are the first low regularity integrators for the derivative nonlinear Schr\"odinger equation. Numerical experiments illustrate our theoretical findings.
\end{abstract}

\maketitle

\section{Introduction}
	
Let us consider the derivative nonlinear Schr\"odinger (dNLS) equation on the one-dimensional torus $\mathbb{T} = \mathbb{R}/2\pi\mathbb{Z}$:

\begin{equation}\label{dnls}
\left\{
\begin{aligned}
&\partial_t u -i \partial_x^2 u = \lambda \partial_x\bigl(|u|^2 u\bigr), \quad t\ge0, \ x \in \mathbb{T},\\
&u(0,x) = u_0(x),
\end{aligned}
\right.
\end{equation}
where $u(t,x):[0,T]\times\mathbb{T}\to\mathbb{C}$ denotes a complex-valued wave envelope. The parameter  $\lambda \in \mathbb{R}$  can be scaled to unity via the transformation $u(t,x) \mapsto |\lambda|^{-1/2}\, u\bigl(t,x \operatorname{sign}\lambda\bigr)$; therefore, we set $\lambda = 1$ without loss of generality.

The dNLS equation was first derived by Mj{\o}lhus~\cite{mjolhus76} as a model for Alfv\'en waves in magnetized plasmas~\cite{mio76,mjolhus89} and has since found applications in nonlinear optics and water-wave theory~\cite{dysthe79,moses07}. Unlike the classical cubic nonlinear Schr\"odinger (NLS) equation, dNLS uses a spatial derivative inside the nonlinearity, which fundamentally alters the analysis. In particular, the derivative causes a genuine loss of regularity in naive energy or fixed-point arguments: estimates that suffice for the semilinear cubic NLS typically fail at the same regularity for dNLS. It also shifts the scaling balance, making low-regularity well-posedness and long-time dynamics more delicate, and it often necessitates the use of structural cancellations (e.g., the gauge transformation \cite{herr06} or the R-transformation \cite{haya93}) together with sharper multilinear estimates.
This change is already visible at the level of conserved quantities, where the additional derivative modifies the energy structure. For instance, although the mass
\begin{equation}\label{mass}
M(u) = \int_{\mathbb{T}} |u|^2 \, dx
\end{equation}
is conserved, the corresponding energy functional for solutions $u\in\mathcal{C}^1([0,T],H^1(\mathbb{T}))$ becomes
\begin{equation}\label{hami}
E(u) = \int_{\mathbb{T}}|u_x|^2+\frac32\,\im\bigl(|u|^2u\overline{u}_x\bigr)+\frac12|u|^6\,dx,
\end{equation}
which involves mixed cubic-quintic interactions and is not sign-definite unless the data are sufficiently small. Moreover, the loss of Galilean invariance and the presence of derivative nonlinearities introduce resonant interactions that are absent in the cubic NLS, posing considerable challenges for both analytical and numerical treatments~\cite{bahouri22,colli01,deng21,li18,nahmod12}.

A powerful tool for handling these difficulties are gauge transformations~\cite{WuYang75,Yang54,Yang74}, which rewrite the equation in a form with a more favorable multilinear structure and partially restored symmetries. On the whole real line $\mathbb{R}$, the gauge transformation introduced by Hayashi and Ozawa~\cite{haya92,haya93} leads to optimal local well‑posedness in $H^s(\mathbb{R})$ for $s\ge\frac12$; see~\cite{tak99}. The transformed equation, although still containing derivative-type terms such as $u^2\overline{u}_x$, becomes amenable to analysis via Bourgain spaces~\cite{bourg93} combined with local smoothing and Strichartz estimates. On the other hand, global well-posedness was first obtained for small data~\cite{colli01,colli02,wu14,wu15} and later extended to arbitrary data in $H^s(\mathbb{R})$ for $s\ge\frac12$ by Bahouri et al.~\cite{bahouri22}. We note that for  $s<\frac12$, the Cauchy problem is ill-posed in the sense that the solution map fails to be uniformly continuous~\cite{biag01}, even though some scaling-critical results exist under additional assumptions; see~\cite{deng21,harrop22}.

The periodic setting presents additional obstacles because dispersion is weaker: local smoothing estimates are unavailable, and Strichartz estimates on the torus exhibit a derivative loss beyond $L^4$ integrability. Consequently, multilinear estimates, in particular, sharp trilinear estimates for terms like $u^2\overline{u}_x$ become delicate in periodic Bourgain spaces. Herr~\cite{herr06} overcame this by proving a sharp trilinear estimate using refined multiplier bounds and a periodic $L^4$ Strichartz estimate. Combined with a periodic gauge transformation, this yields local well-posedness in $H^s(\mathbb{T})$ for $s\ge\frac12$ and, under a small-$L^2$ condition, global well-posedness for $s\ge1$; see~\cite{bahouri24}.

In contrast, numerical studies of \eqref{dnls} are comparatively scarce. The derivative nonlinearity makes explicit schemes particularly challenging. Although implicit finite-difference methods have been proposed by Li et al.~\cite{li17,li18}, their convergence relies on the assumption of sufficiently smooth solutions. Designing efficient explicit schemes that remain effective for low regularity data therefore remains an open problem.

Recently, low regularity error analysis for dispersive equations has attracted considerable interest, with developments for the NLS~\cite{bao24,baomw,jisiam,liwunls,ostfocm,ostjems}, the KdV equation~\cite{caolwy,hof17,liwzgauge,liwukdv,maiermcom,wuintpart}, and other dispersive models~\cite{caolly,jimcom,jilisu,jizhao,swz21,wangz22}. However, low regularity integrators and their rigorous error analysis for the dNLS equation are still lacking.

In this work, we will construct a first-order exponential-type integrator for the periodic dNLS equation \eqref{dnls}. Unlike traditional methods~\cite{lubich08}, we will prove first-order convergence in $H^s(\mathbb{T})$ for initial data in $H^{s+1}$. Moreover, we will construct a second-order symmetric method that also exhibits first-order convergence for low regularity but has a better long-term behavior. The analysis is carried out entirely within a Sobolev framework and no filtering techniques or CFL-type conditions are required.

\subsection*{Outline of the paper.} The rest of this paper is organized as follows. In Section~\ref{sectiongauge}, we recall the gauge transformation introduced in \cite{herr06} for the dNLS~\eqref{dnls} on the torus, which enables us to solve the Cauchy problem and derive low regularity integrators for the dNLS. In Section~\ref{sectionmethod}, we introduce two exponential-type integrators and give the main convergence result (Theorem~\ref{mainthm}). The local and global error estimates for the methods are given in Sections~\ref{sectionlocal} and \ref{sectionglobal}, respectively. Numerical verification of theoretical results are provided at the end of the paper.
\newpage
\subsection*{Notations.}
\begin{itemize}[leftmargin=*]
\item For $a, b\ge 0$, we write \( a \lesssim b \) to indicate that there exists a constant \( C > 0 \), independent of the time step size \( \tau \in (0,1) \), such that \( a \le Cb \). If the constant depends in particular on a parameter \( c \), we write \( \lesssim_c \). The notation \( a \sim b \) means both \( a \lesssim b \) and \( b \lesssim a \).
		
\item For $y\in\mathbb{R}$, we use the Japanese bracket $\langle y\rangle=(1+|y|^2)^{\frac12}$.
		
\item For any $s\in\mathbb{R}$ and a function $f(x)=\sum\limits_{k\in\mathbb{Z}}\widehat{f}_k\fe^{ikx}$, we define the Sobolev norm for $H^s(\mathbb{T})$ as 
$$
\|f\|_s^2=2\pi\sum_{k\in\mathbb{Z}}(1+k^2)^s|\widehat{f}_k|^2.
$$
For the sake of brevity, we denote by $H^s$ the spatial functions defined on $\mathbb{T}$ equipped with the above norm $\|\cdot\|_s$. Moreover, for $s=0$, the space reduces to $L^2$ and the corresponding norm, denoted as $\|\cdot\|_0$, is derived from the inner product
$$
\langle f, g \rangle =2\pi \sum_{k \in \mathbb{Z}} \widehat{f}_k \overline{\widehat{g}}_k = \int_0^{2\pi} f(x) \overline{g}(x) \, dx.
$$
\item For $u \in C([-T,T],L^2(\mathbb{T}))$, we define 
$$
\Pi_0u=\frac{1}{2\pi}\int_0^{2\pi}u\,dx
$$
and set $\mu(u)=\Pi_0|u|^2$. Due to mass conservation~\eqref{mass}, $\mu(\cdot)$ is constant along trajectories. Therefore, we write $\mu=\mu(u_0)$ for short, where $u_0$ is the initial data of the trajectory.
		
\item For $c\in \mathbb{R}$, we define the translations $\kappa_{c}u(t,x):=u(t,x+2ct)$ for $u\in C([-T,T],L^2(\mathbb{T}))$.
\end{itemize}
	
\section{Gauge Transformation and Cauchy Problem for \eqref{dnls}}\label{sectiongauge}

In this section, we shall recall a gauge transformation for the dNLS~\eqref{dnls} on the torus and formulate the Cauchy problem for \eqref{dnls}. This transformation allows us to transfer the derivative from the entire nonlinearity to the conjugate term, $\overline{u}$. This is a crucial step in both analyzing the Cauchy problem and deriving low regularity integrators. For more details of this transformation and the Cauchy problem, we refer to~\cite{herr06}.

\subsection{The gauge transformation}

We first recall the gauge transformation from~\cite{herr06}.

\begin{definition}\label{defgauge}
For $f \in L^2(\mathbb{T})$, we set
$$
\mathcal{G}_0(f)(x)=\fe^{-i\mathcal{I}(f)}f(x),
$$
where the phase \(\mathcal{I}(f)\) is given by
$$
\mathcal{I}(f)(x)=\frac{1}{2\pi}\int_0^{2\pi}\int_{\xi}^x |f(y)|^2
-\frac{1}{2\pi}\|f\|_0^2 \, dy\, d\xi.
$$
Then, for $u\in C([0,T],L^2(\mathbb{T}))$, we define
\begin{equation}\label{eqgaugedef}
\mathcal{G}(u)(t,x)=\kappa_{-\mu(u)}\mathcal{G}_0(u(t))(x).
\end{equation}
\end{definition}
Since $|f(y)|^2-\frac{1}{2\pi}\|f\|_0^2$ is a zero-mean function on $\mathbb{T}$, its primitive $\mathcal{I}(f)(x)$ is $2\pi$-periodic. As a consequence, the transformed function $\mathcal{G}_0(f)(x) = \fe^{-i\mathcal{I}(f)(x)}f(x)$ is also $2\pi$-periodic.

\begin{lemma}[{\cite[Lemma 2.1]{herr06}}]\label{lemhomo}
For any $\gamma\ge0$, the transformation
$$
\mathcal{G}:\mathcal{C}([0,T],H^\gamma(\mathbb{T}))\to\mathcal{C}([0,T],H^\gamma(\mathbb{T}))
$$
is a gauge transformation, i.e., the map $\mathcal{G}$ is a homeomorphism that satisfies
\begin{equation}\label{gaugeprop}
\|\mathcal{G}(u)(t,\cdot)-\mathcal{G}(v)(t,\cdot)\|_\gamma\leq C\|u(t,\cdot)-v(t,\cdot)\|_\gamma,\quad t \in [0,T],
\end{equation}
where $C$ depends on $\|u\|_{L^\infty([0,T],H^\gamma(\mathbb{T}))}$ and $\|v\|_{L^\infty([0,T], H^\gamma(\mathbb{T}))}$.
Furthermore, the inverse map is given by
$$
\mathcal{G}^{-1}(v)(t,x)=\fe^{i\mathcal{I}(\kappa_{\mu(v)}v)}\kappa_{\mu(v)} v(t,x)
$$
and $\mathcal{G}^{-1}$ satisfies a similar estimate as \eqref{gaugeprop} on $[0,T]$.
\end{lemma}

Note that $\mu(u)=\mu(v)$ for $v=\mathcal{G}(u)$. Moreover, with the help of this gauge transformation, we have the following result:

\begin{corollary}\label{corogauge}
Let $u,v\in\mathcal{C}([0,T],H^\gamma(\mathbb{T}))$ with $\gamma\ge\frac12$ and $v=\mathcal{G}(u)$.
Then, $u$ is a mild solution of \eqref{dnls} if and only if $v$ is a mild solution of
\begin{equation}\label{gdnls}
\left\{\begin{aligned}
&\partial_t v-i\partial^2_{x}v=F(v) \quad\textrm{in }[0,T]\times\mathbb{T},\\
&v(0)=\mathcal{G}(u_0)=\mathcal{G}_0(u_0)=v_0,
\end{aligned}\right.
\end{equation}
where
\begin{equation}\label{Fpsi}
\begin{aligned}
F(v)&=-v^2\overline{v}_x+\tfrac{i}{2}|v|^4v-i\mu|v|^2v+i\psi(v)v,\\
\psi(v)&=\Pi_0\bigl(2\im(\overline{v}_x v)(t)-\frac12|v|^4(t)\bigr)+\mu^2 \ \text{with} \ \mu=\frac{1}{2\pi}\|v_0\|_0^2.
\end{aligned}
\end{equation}
\end{corollary}

In the following, we denote
\begin{alignat}{2}
F_1(v) &= -v^2\overline{v}_x,\qquad\qquad\qquad &&F_2(v) = \frac i2|v|^4v,\nonumber\\
\label{f1234}F_3(v) &= -i\mu|v|^2v,&&F_4(v) = i\psi(v)v.
\end{alignat}
Thus, Duhamel's formula for the transformed equation \eqref{gdnls} reads
\begin{equation}\label{duhv}
v(t)=\fe^{it\partial_x^2}v_0+\int_0^t\sum\limits_{j=1}^4\fe^{i(t-\nu)\partial_x^2}F_j\bigl(v(\nu)\bigr)\,d\nu.
\end{equation}

\subsection{Cauchy problem for (\ref{dnls})}

Thanks to Corollary~\ref{corogauge}, we can obtain local well-posedness of \eqref{dnls} for $u_0\in H^\gamma$ with $\gamma\ge\frac12$. Indeed, let us define the Bourgain space \cite{bourg93,herr06} by its norm
$$
\|u\|_{Z^\gamma(\mathbb{R}\times\mathbb{T})}=\|\langle k\rangle^\gamma\langle\sigma+k^2\rangle^\frac12\widetilde{u}(\sigma,k)\|_{L^2(\sigma)l^2(k)}+\|\langle k\rangle^\gamma\widetilde{u}(\sigma,k)\|_{L^1(\sigma)l^2(k)},
$$
where $\widetilde{u}$ denotes the ``time-space" Fourier transform of $u$. Then, we have the following proposition:

\begin{proposition}[{\cite[Theorem 5.1]{herr06}}]\label{propexist}
For any $u_0\in H^\gamma(\mathbb{T})$ with $\gamma\ge\frac12$, there exists a unique solution $v(t,x)\in Z^\gamma\subset\mathcal{C}([0, T], H^\gamma)$ satisfying the equation~\eqref{gdnls} on $[0,T]$, and $u(t,x)=\mathcal{G}^{-1}(v)(t,x)$ is the unique continuous (in time) solution to \eqref{dnls} in the complete metric space
$$
Y^\gamma=\{u=\mathcal{G}(v):v\in Z^\gamma\}.
$$
Moreover, we have
\begin{equation}\label{uvbound}
\sup\limits_{t\in[0,T]}\|u(t,\cdot)\|_\gamma\le C_1\sup\limits_{t\in[0,T]}\|v(t,\cdot)\|_\gamma\le C_2\|v\|_{L^\infty H^\gamma}\le C_3\|v(t,\cdot)\|_{Z^\gamma}\leq C_T,
\end{equation}
where $C_T$ depends on $\gamma$, $T$ and $\|u_0\|_\gamma$.
\end{proposition}

Note that here $u$ and $v$ are defined globally, but are solutions to \eqref{dnls} and \eqref{gdnls}, respectively, only on $[0,T]$~\cite{herr06,ostjems}. For the proof of the proposition, we refer to \cite{herr06}. Note that the estimate~\eqref{uvbound} follows from the embedding $L^1(\sigma)H^\gamma\subset L^\infty(t)H^\gamma$, i.e.,
$$
\sup_{t\in[0,T]}\|f(t)\|_{H^\gamma}\lesssim\inf_{g\mid_{[0,T]}=f\mid_{[0,T]}}\int_{\mathbb R}\|\widehat{g}(\sigma)\|_{H^\gamma}\,d\sigma,
$$
where $\widehat{f}(\sigma)$ denotes the Fourier transform of $f(t)$ with respect to $t$.

\section{Construction of the Numerical Method}\label{sectionmethod}

In this section, we construct numerical methods for the dNLS equation~\eqref{dnls} and state the main result of the paper, Theorem~\ref{mainthm}. We first introduce a basic first-order method and then extend it to a symmetric method that performs better in practice.

\subsection{Construction of the basic method}\label{sectmeth1}

First, we derive a basic low-regularity exponential integrator. Inspired by \cite{ostfocm}, we employ the twist $z=\fe^{-it\partial_x^2}v$, where $v$ is defined in Proposition~\ref{propexist}. Using Duhamel's formula \eqref{duhv}, we get
\begin{equation}\label{duhz}
z(t_{n+1})=z(t_n)+\int_0^\tau\sum\limits_{j=1}^4\fe^{-i(t_n+\nu)\partial_x^2}F_j\Bigl(\fe^{i(t_n+\nu)\partial_x^2}z(t_n+\nu)\Bigr)\,d\nu.
\end{equation}

We first consider the term containing $F_1$. We approximate $z(t_n+\nu)\approx z(t_n)$ in \eqref{duhz} and use the Fourier expansion 
$$
z(t)=\sum\limits_{k\in\mathbb{Z}}\widehat{z}_k(t)\fe^{ikx}
$$
to get
\begin{equation}\label{kk1}
\begin{aligned}
\int_0^\tau&\fe^{-i(t_n+\nu)\partial_x^2}F_1\bigl(\fe^{i(t_n+\nu)\partial_x^2}z(t_n+\nu)\bigr)\,d\nu\approx\int_0^\tau\fe^{-i(t_n+\nu)\partial_x^2}F_1\bigl(\fe^{i(t_n+\nu)\partial_x^2}z(t_n)\bigr)\,d\nu\\
&\qquad=i\sum\limits_{k\in\mathbb{Z}}\sum\limits_{\substack{k_1,k_2,k_3\in\mathbb{Z}\\-k_1+k_2+k_3=k}}\int_0^\tau k_1\fe^{i(t_n+\nu)(k^2+k_1^2-k_2^2-k_3^2)}\overline{\widehat{z}}_{k_1}(t_n)\widehat{z}_{k_2}(t_n)\widehat{z}_{k_3}(t_n)\,d\nu\,\fe^{ikx}.
\end{aligned}
\end{equation}
For getting a formula in physical space it suffices to estimate the exponential kernel,
$$
\int_0^\tau\fe^{i\nu(k^2+k_1^2-k_2^2-k_3^2)}\,d\nu,
$$
in an appropriate way and to transform back.
Note that
$$
k^2+k_1^2-k_2^2-k_3^2=-2kk_1+2k_2k_3.
$$
To capture the dominant derivative contribution in the oscillatory phase, we use the approximation $\fe^{i\nu(k^2+k_1^2-k_2^2-k_3^2)}\approx \fe^{-2i\nu kk_1}$. Since the sum $\eqref{kk1}$ vanishes for $k_1=0$, we thus have
\begin{align*}
\sum\limits_{k\in\mathbb{Z}}&\sum\limits_{\substack{k_1,k_2,k_3\in\mathbb{Z}\\-k_1+k_2+k_3=k}}\int_0^\tau k_1\fe^{i(t_n+\nu)(k^2+k_1^2-k_2^2-k_3^2)}\overline{\widehat{z}}_{k_1}(t_n)\widehat{z}_{k_2}(t_n)\widehat{z}_{k_3}(t_n)\,d\nu\,\fe^{ikx}\\
&\approx\sum\limits_{k\in\mathbb{Z}\setminus\{0\}}\sum\limits_{\substack{k_1,k_2,k_3\in\mathbb{Z}\\-k_1+k_2+k_3=k}}\fe^{it_n(k^2+k_1^2-k_2^2-k_3^2)}\frac{\fe^{-2i\tau kk_1}-1}{-2ik}\overline{\widehat{z}}_{k_1}(t_n)\widehat{z}_{k_2}(t_n)\widehat{z}_{k_3}(t_n)\fe^{ikx}\\
&\quad+\tau\sum\limits_{\substack{k_1,k_2,k_3\in\mathbb{Z}\\k_1=k_2+k_3}}k_1\fe^{it_n(k_1^2-k_2^2-k_3^2)}\overline{\widehat{z}}_{k_1}(t_n)\widehat{z}_{k_2}(t_n)\widehat{z}_{k_3}(t_n).
\end{align*}
Twisting back and using $-2kk_1=k^2+k_1^2-(k_2+k_3)^2$, we obtain
\begin{equation}\label{getg1}
\begin{aligned}
\int_0^\tau&\fe^{-i(t_n+\nu)\partial_x^2}F_1\bigl(\fe^{i(t_n+\nu)\partial_x^2}z(t_n+\nu)\bigr)\,d\nu\\
&\approx\fe^{-it_n\partial_x^2}\sum\limits_{k\in\mathbb{Z}\setminus\{0\}}\sum\limits_{\substack{k_1,k_2,k_3\in\mathbb{Z}\\-k_1+k_2+k_3=k}}\frac{\fe^{i\tau(k^2+k_1^2-(k_2+k_3)^2)}-1}{-2k}\overline{\widehat{v}}_{k_1}(t_n)\widehat{v}_{k_2}(t_n)\widehat{v}_{k_3}(t_n)\fe^{ikx}\\
&\quad+i\tau\sum\limits_{\substack{k_1,k_2,k_3\in\mathbb{Z}\\k_1=k_2+k_3}}k_1\overline{\widehat{v}}_{k_1}(t_n)\widehat{v}_{k_2}(t_n)\widehat{v}_{k_3}(t_n)\\
&=-\frac i2\partial_x^{-1}\fe^{-it_n\partial_x^2}g_0\bigl(v(t_n),v(t_n),\tau\bigr)-\tau\Pi_0(v^2(t_n)\overline{v}_x(t_n)),
\end{aligned}
\end{equation}
where
\begin{subequations}
\begin{align}
\label{g0diff} g_0(v_1,v_2,\tau)&=(I-\Pi_0)\Bigl(\fe^{-i\tau\partial_x^2}\bigl((\fe^{-i\tau\partial_x^2}\overline{v}_1)\fe^{i\tau\partial_x^2}v_2^2\bigr)-\overline{v}_1v_2^2\Bigr)\\
\label{g0int} &=-2i\partial_x\int_0^\tau\fe^{-i\nu\partial_x^2}\bigl((\partial_x\fe^{-i\nu\partial_x^2}\overline{v}_1)\fe^{i\nu\partial_x^2}v_2^2\bigr)\,d\nu.
\end{align}
\end{subequations}
In the following, we will also denote
\begin{equation}\label{g1ofu}
g_1(v,\tau)=g_0(v,v,\tau)
\end{equation}
for short.

The approximation \eqref{getg1}, however, is not sufficient for $F_1$ at low regularity. As a compensation, we approximate the term containing $F_2$ as follows:
\begin{align*}
\int_0^\tau&\fe^{-i(t_n+\nu)\partial_x^2}F_2\bigl(\fe^{i(t_n+\nu)\partial_x^2}z(t_n+\nu)\bigr)\,d\nu\\
&\qquad\approx\frac i2\tau\fe^{-it_n\partial_x^2}\Bigl((\fe^{-it_n\partial_x^2}z(t_n))^2\varphi_1(-2i\tau\partial_x^2)\bigl((\fe^{it_n\partial_x^2}\overline{z}(t_n))^2\fe^{-it_n\partial_x^2}z(t_n)\bigr)\Bigr),
\end{align*}
where
\begin{equation}\label{defphi1}
\varphi_1(z)=\frac{e^z-1}{z}
\end{equation}
as defined in \cite{hoch05,hochacta}. Note that this approximation is different from that given in \cite{ostfocm}. The reason for this will become clear in the proof of Proposition~\ref{proplocal1}.

For $F_3$, we take the approximation given in \cite{ostfocm}
$$
\int_0^\tau\fe^{-i(t_n+\nu)\partial_x^2}F_3\bigl(\fe^{i(t_n+\nu)\partial_x^2}z(t_n+\nu)\bigr)\,d\nu\approx-i\mu\tau\fe^{-it_n\partial_x^2}\bigl((\fe^{-it_n\partial_x^2}z(t_n))^2\varphi_1(-2i\tau\partial_x^2)\fe^{it_n\partial_x^2}\overline{z}(t_n)\bigr).
$$

For $F_4$, since $\psi$ is independent of $x$ and
$$
\psi(v)=\psi(\fe^{it\partial_x^2}z)=\psi(z),
$$
we can use
$$
\int_0^\tau\fe^{-i(t_n+\nu)\partial_x^2}F_4\bigl(\fe^{i(t_n+\nu)\partial_x^2}z(t_n+\nu)\bigr)\,d\nu\approx i\tau\psi(z(t_n))z(t_n).
$$

Combining the above approximations of all $F_i,~i=1,2,3,4$, we obtain a first-order method for ~\eqref{gdnls}:
\begin{equation}\label{meth1}
\begin{aligned}
v_{n+1}=\Psi_1^\tau(v_n)&=\fe^{i\tau\partial_x^2}\Bigl(v_n-\frac i2\partial_x^{-1}g_1(v_n,\tau)-\tau\Pi_0(v_n^2\partial_x\overline{v}_n)\\
&\quad+\frac i2\tau v_n^2\varphi_1(-2i\tau\partial_x^2)(\overline{v}_n^2v_n)-i\mu\tau v_n^2\varphi_1(-2i\tau\partial_x^2)\overline{v}_n+i\tau\psi(v_n)v_n\Bigr),
\end{aligned}
\end{equation}
where $\psi,g_1,\varphi_1$ are defined in \eqref{Fpsi}, \eqref{g1ofu} and \eqref{defphi1}, respectively.

\subsection{Construction of a symmetric method}\label{sectmeth2}

Due to the use of the gauge transformation~\eqref{eqgaugedef}, the conservation behavior, especially of mass, plays a significant role in the numerical solution of the dNLS equation. To obtain a better conservation behavior, we will construct a symmetric method for the transformed system~\eqref{gdnls}.

Since the equation is time-reversible, we can consider a two-step method that treats $v_{n+1}$ and $v_{n-1}$ symmetrically. Replacing $\tau$ by $-\tau$ in the basic method \eqref{meth1}, we get
\begin{align*}
v_{n-1}=\Psi_1^{-\tau}(v_n)&=\fe^{-i\tau\partial_x^2}\Bigl(v_n-\frac i2\partial_x^{-1}g_1(v_n,-\tau)+\tau\Pi_0(v_n^2\partial_x\overline{v}_n)\\
&\quad-\frac i2\tau v_n^2\varphi_1(2i\tau\partial_x^2)(\overline{v}_n^2v_n)+i\mu\tau v_n^2\varphi_1(2i\tau\partial_x^2)\overline{v}_n-i\tau\psi(v_n)v_n\Bigr).
\end{align*}
Combining this with \eqref{meth1}, we obtain a second-order symmetric method \cite[Section~XV.1]{hairer}
\begin{subequations}\label{meth2}
\begin{equation}
\begin{aligned}
v_{n+1}&=\Psi_2^\tau(v_n,v_{n-1})=\fe^{2i\tau\partial_x^2}v_{n-1}+\fe^{i\tau\partial_x^2}\Bigl(-\frac{i}2\partial_x^{-1}g_2(v_n,\tau)-2\tau\Pi_0(v_n^2\partial_x\overline{v}_n)\\
&\quad+i\tau v_n^2\fe^{2i\tau\partial_x^2}\varphi_1(-4i\tau\partial_x^2)(\overline{v}_n^2v_n)
-2i\mu\tau v_n^2\fe^{2i\tau\partial_x^2}\varphi_1(-4i\tau\partial_x^2)\overline{v}_n+2i\tau\psi(v_n)v_n\Bigr),
\end{aligned}
\end{equation}
where
$$
 g_2(v,\tau)=g_1(v,\tau)-g_1(v,-\tau).
$$
It remains to define $v_1$. In order to keep the symmetry of the method, we use the implicit scheme
\begin{equation}
v_1=(\Psi_{1}^{-\frac{\tau}{2}})^{-1}\circ\Psi_1^{\frac{\tau}{2}}(v_0).
\end{equation}
\end{subequations}
Note that the computational cost of this method is almost the same as that of the basic method~\eqref{meth1}.

\subsection{Main theorem}

We are now ready to present the main theorem of this paper. We will only present the convergence result for the basic method~\eqref{meth1}. Despite being a second-order method, the symmetric method~\eqref{meth2} only converges with first order under the regularity constraints of Theorem~\ref{mainthm}. We will omit the details of the error analysis for ~\eqref{meth2} since the analysis is very similar to that of \eqref{meth1}. Moreover, in Section~\ref{sectionnumexp}, we will demonstrate that the symmetric method performs better than the basic method in practice.

\begin{theorem}\label{mainthm}
For given $s>\frac12$ and initial data $u_0\in H^{s+1}(\mathbb{T})$, let $u(t,x)$ be the exact solution to~\eqref{dnls} given by Corollary~\ref{corogauge} and Proposition~\ref{propexist}. Further, let $v_n$ be the sequence defined in \eqref{meth1} with initial data $v_0=\mathcal{G}_0(u_0)$. Then, there exist constants $\tau_0,~C_T>0$ such that, for any $\tau\in(0,\tau_0]$, we have the error estimate
\begin{equation}\label{estorder1}
\|\mathcal{G}^{-1}(v_n)-u(t_n)\|_s\leq C_T\tau,\qquad 0\leq t_n=n\tau\leq T,
\end{equation}
where $\mathcal{G}$ is given by Definition~\ref{defgauge}. The constants $\tau_0$ and $C_T$ depend on $\|u_0\|_{s+1}$ and $T$, but they are independent of $\tau$ and $n$.
\end{theorem}
The proof of this theorem is given in Sections~\ref{sectionlocal} and \ref{sectionglobal}.

\section{Local error analysis}\label{sectionlocal}

In this section, we shall start analyzing the first-order method~\eqref{meth1}. We will first present several useful lemmas, then estimate the local error.

\subsection{Some technical tools}

The most useful tool in this work is the well-known classical bilinear estimate
\begin{equation}\label{cbe}
\|fg\|_s\le C_s\|f\|_s\|g\|_s, \qquad s>\frac12.
\end{equation}
Moreover, to better handle the derivative in $F_1$, we shall also use another bilinear estimate stated in \cite[Lemma~2(i)]{liwunls}:
\begin{equation}\label{cbepm}
\|fg\|_{s-1}\le C_s\|f\|_{s-1}\|g\|_s, \qquad s>\frac12.
\end{equation}

The following lemma will be used several times throughout the rest of this paper.
\begin{lemma}\label{lempi0}
For any $f,g\in H^\frac12$, we have
\begin{equation}\label{pi0fgx}
|\Pi_0(f\overline{g}_x)|\le\|f\|_\frac12\|g\|_\frac12.
\end{equation}
\end{lemma}
\begin{proof}
Note that $\Pi_0(f\overline{g}_x)$ is the $L^2$-inner product of $f$ and $g_x$. We thus have
$$
|\Pi_0(f\overline{g}_x)|=|\langle f,\partial_xg\rangle|=|\langle\partial_x^\frac12f,\partial_x\partial_x^{-\frac12}g\rangle|\le\|\partial_x^\frac12f\|_0\|\partial_x^\frac12g\|_0\leq\|f\|_\frac12\|g\|_\frac12,
$$
where we used the self-adjointness of $\partial_x^\frac12$ and the Cauchy--Schwarz inequality in $L^2$. This concludes the proof.
\end{proof}

Next, in order to reduce the derivative loss caused by the approximation $z(t_n+\nu)\approx z(t_n)$, we shall need another lemma, which has been used similarly for the KdV equation in \cite{wuintpart}:
\begin{lemma}[Integration by parts]\label{lemintpart}
For any functions $f$ and $g$ defined on $\mathbb{R}\times\mathbb{T}$, we have the following identity:
\begin{equation}\label{intpart}
\begin{aligned}
(I&-\Pi_0)\int_0^\tau\fe^{-i\nu\partial_x^2}\Bigl(\bigl(\fe^{-i\nu\partial_x^2}\partial_x\overline{f}(\nu)\bigr)\fe^{i\nu\partial_x^2}g(\nu)\Bigr)\,d\nu=\frac{i}2\fe^{-i\tau\partial_x^2}\partial_x^{-1}(I-\Pi_0)\Bigl(\bigl(\fe^{-i\tau\partial_x^2}\overline{f}(\tau)\bigr)\fe^{i\tau\partial_x^2}g(\tau)\Bigr)\\
&-\frac{i}2\partial_x^{-1}(I-\Pi_0)\bigl(\overline{f}(0)g(0)\bigr)-\frac{i}2\int_0^\tau\fe^{-i\nu\partial_x^2}\partial_x^{-1}(I-\Pi_0)\Bigl(\bigl(\fe^{-i\nu\partial_x^2}\partial_\nu\overline{f}(\nu)\bigr)\fe^{i\nu\partial_x^2}g(\nu)\Bigr)\,d\nu\\
&-\frac{i}2\int_0^\tau\fe^{-i\nu\partial_x^2}\partial_x^{-1}(I-\Pi_0)\Bigl(\bigl(\fe^{-i\nu\partial_x^2}\overline{f}(\nu)\bigr)\fe^{i\nu\partial_x^2}\partial_\nu g(\nu)\Bigr)\,d\nu.
\end{aligned}
\end{equation}
\end{lemma}
\begin{proof}
Taking the spatial Fourier transform, we obtain
\begin{align*}
\mathcal{F}_x\biggl((I&-\Pi_0)\int_0^\tau\fe^{-i\nu\partial_x^2}\Bigl(\bigl(\fe^{-i\nu\partial_x^2}\partial_x\overline{f}(\nu)\bigr)\fe^{i\nu\partial_x^2}g(\nu)\Bigr)\,d\nu\biggr)_k\\
&=\begin{cases}\displaystyle
-i\int_0^\tau\sum_{\substack{k_1,k_2\in\mathbb{Z}\\-k_1+k_2=k}}k_1\fe^{i\nu\alpha}\overline{\widehat{f}}(\nu,k_1)\widehat{g}(\nu,k_2)\,d\nu, & k\neq0,\\
0, & k=0,
\end{cases}
\end{align*}
where
$$
\alpha=k^2+k_1^2-k_2^2=-2kk_1.
$$
Integrating by parts in time, we obtain
\begin{align*}
-i\int_0^\tau&\sum_{\substack{k_1,k_2\in\mathbb{Z}\\-k_1+k_2=k\neq0}}k_1\fe^{i\nu\alpha}\overline{\widehat{f}}(\nu,k_1)\widehat{g}(\nu,k_2)\,d\nu\\
&=\sum_{\substack{k_1,k_2\in\mathbb{Z}\\-k_1+k_2=k\neq0}}\frac1{2k}\fe^{i\nu\alpha}\overline{\widehat{f}}(\nu,k_1)\widehat{g}(\nu,k_2)\Big|_0^\tau-\int_0^\tau\sum_{\substack{k_1,k_2\in\mathbb{Z}\\-k_1+k_2=k\neq0}}\frac1{2k}\fe^{i\nu\alpha}\partial_\nu\bigl(\overline{\widehat{f}}(\nu,k_1)\widehat{g}(\nu,k_2)\bigr)\,d\nu.
\end{align*}
Note that the right-hand side of the above equation vanishes when $k_1=0$. We then conclude by multiplying with $\fe^{ikx}$ and summing over $k\neq0$.
\end{proof}

\subsection{Local error analysis}

We now step to the local error. By \eqref{meth1} and Duhamel's formula, we first write the local error as follows:
\begin{equation}\label{vloc1}
\begin{aligned}
v&(t_{n+1})-\Psi_1^\tau(v(t_n))\\
&=\int_0^\tau\fe^{i(\tau-\nu)\partial_x^2}F_1(v(t_n+\nu))\,d\nu+\fe^{i\tau\partial_x^2}\Bigl(\frac{i}2\partial_x^{-1}g_1(v(t_n),\tau)+\tau\Pi_0\bigl(v^2(t_n)\overline{v}_x(t_n)\bigr)\Bigr)\\
&\quad+\int_0^\tau\fe^{i(\tau-\nu)\partial_x^2}F_2(v(t_n+\nu))\,d\nu-\frac{i}2\tau\fe^{i\tau\partial_x^2}\Bigl(v^2(t_n)\varphi_1(-2i\tau\partial_x^2)\bigl(\overline{v}^2(t_n)v(t_n)\bigr)\Bigr)\\
&\quad+\int_0^\tau\fe^{i(\tau-\nu)\partial_x^2}F_3(v(t_n+\nu))\,d\nu +i\mu\tau\fe^{i\tau\partial_x^2}\Bigl(v^2(t_n)\varphi_1(-2i\tau\partial_x^2)\overline{v}(t_n)\Bigr)\\
&\quad+\int_0^\tau\fe^{i(\tau-\nu)\partial_x^2}F_4(v(t_n+\nu))\,d\nu-i\tau\psi(v(t_n))\fe^{i\tau\partial_x^2}v(t_n)\\
&=\mathcal{E}_1(t_n)+\mathcal{E}_2(t_n)+\mathcal{E}_3(t_n)+\mathcal{E}_4(t_n),
\end{aligned}
\end{equation}
where $F_j,~j=1,2,3,4$ and $g_1$ are defined in \eqref{f1234} and \eqref{g1ofu}, respectively.

Next, we shall estimate the error caused by the approximation $z(t_n+\nu)\approx z(t_n)$, and prove the following lemma.
\begin{lemma}\label{lemvdiff}
For any $s>\frac12$ and $v(t,x)\in Z^{s+1}\hookrightarrow\mathcal{C}([0,T],H^{s+1})$ as in Proposition~\ref{propexist}, we have for any $\nu\in(0,\tau]$ with $\tau$ sufficiently small that
\begin{equation}\label{vdiff}
\|v(t_n+\nu)-\fe^{i\nu\partial_x^2}v(t_n)\|_s\le C_T\nu,\qquad 0\le t_n=n\tau\le T,
\end{equation}
where $C_T$ depends on $s$ and $\sup\limits_{\xi\in[0,\nu]}\|v(t_n+\xi)\|_{s+1}$.
\end{lemma}
\begin{proof}
By Duhamel's formula \eqref{duhv}, it suffices to prove
$$
\sup\limits_{\xi\in[0,\nu]}\|F_j(v(t_n+\xi))\|_s\le C_T,\qquad j=1,2,3,4
$$
since $\fe^{it\partial_x^2}$ is an isometry in $H^s$. Indeed, by \eqref{uvbound} and \eqref{cbe}, we have
\begin{align*}
\sup\limits_{\xi\in[0,\nu]}\|F_1(v(t_n+\xi))\|_s&\lesssim\sup\limits_{\xi\in[0,\nu]}\|v(t_n+\xi)\|^2_s\|v(t_n+\xi)\|_{s+1}\lesssim1,\\
\sup\limits_{\xi\in[0,\nu]}\|F_2(v(t_n+\xi))\|_s&\lesssim\sup\limits_{\xi\in[0,\nu]}\|v(t_n+\xi)\|^5_s\lesssim1.
\end{align*}
Moreover, noting that $2\pi\mu=\|u_0\|_0^2=\|v_0\|_0^2\le\|v_0\|^2_s$, we obtain
$$
\sup\limits_{\xi\in[0,\nu]}\|F_3(v(t_n+\xi))\|_s\lesssim\sup\limits_{\xi\in[0,\nu]}\|v_0\|^2_s\|v(t_n+\xi)\|^3_s\lesssim1.
$$
For $F_4$, by using Lemma~\ref{lempi0}, we get from Sobolev embedding
$$
\sup\limits_{\xi\in[0,\nu]}|\psi(v(t_n+\xi))|\lesssim\sup\limits_{\xi\in[0,\nu]}\bigl(\|v(t_n+\xi)\|^2_\frac12+\|v(t_n+\xi)\|^4_{L^4}+\|v_0\|^2_s\bigr)\lesssim1+\sup\limits_{\xi\in[0,\nu]}\|v(t_n+\xi)\|^4_\frac14\lesssim1.
$$
This yields
$$
\sup\limits_{\xi\in[0,\nu]}\|F_4(v(t_n+\xi))\|_s\lesssim\sup\limits_{\xi\in[0,\nu]}|\psi(v(t_n+\xi))|\|v(t_n+\xi)\|_s\lesssim1,
$$
and concludes the proof.
\end{proof}

We now state the local error result.

\begin{proposition}\label{proplocal1}
For any $s>\frac12$ and $v(t,x)\in Z^{s+1}$ as in Proposition~\ref{propexist}, we have for $\tau$ sufficiently small the following bound for the local error:
\begin{equation}\label{estloc1}
\|v(t_{n+1})-\Psi_1^\tau(v(t_n))\|_s\le C_T\tau^2,\qquad 0\le t_n=n\tau\le T,
\end{equation}
where $C_T$ depends on $s$ and $\sup\limits_{\nu\in[0,\tau]}\|v(t_n+\nu)\|_{s+1}$.
\end{proposition}
\begin{proof}
By \eqref{vloc1}, it suffices to prove that
$$
\|\mathcal{E}_j(t_n)\|_s\lesssim\tau^2,\qquad j=1,2,3,4.
$$

In the following, we denote
\begin{align}\label{wnuandv}
w(\nu)=\fe^{-i\nu\partial_x^2}v(t_n+\nu),\qquad \nu\in[0,\tau]
\end{align}
with its Fourier expansion
$$
w(\nu)=\sum\limits_{k\in\mathbb{Z}}\widehat{w}_k(\nu)\fe^{ikx}.
$$
Thus, by Lemma~\ref{lemvdiff}, we directly obtain
\begin{equation}\label{wdiff}
\sup_{\nu\in[0,\tau]}\|w(\nu)-w(0)\|_s\lesssim\tau.
\end{equation}

We first estimate $\mathcal{E}_1$. From the calculations in Sect.~\ref{sectmeth1}, we infer that
\begin{align*}
-\frac{i}2\fe^{i\tau\partial_x^2}&\partial_x^{-1}g_1(v(t_n),\tau)-\fe^{i\tau\partial_x^2}\tau\Pi_0(v^2(t_n)\overline{v}_x(t_n))\\
&=i\fe^{i\tau\partial_x^2}\int_0^\tau\sum\limits_{k\in\mathbb{Z}}\sum\limits_{\substack{k_1,k_2,k_3\in\mathbb{Z}\\-k_1+k_2+k_3=k}}\fe^{-2i\nu kk_1}k_1\overline{\widehat{w}}_{k_1}(0)\widehat{w}_{k_2}(0)\widehat{w}_{k_3}(0)\,d\nu\,\fe^{ikx}.
\end{align*}
Therefore, we rewrite $\mathcal{E}_1(t_n)$ in its Fourier expansion
\begin{align*}
\mathcal{E}_1(t_n)&=i\fe^{i\tau\partial_x^2}\int_0^\tau\sum\limits_{k\in\mathbb{Z}}\sum\limits_{\substack{k_1,k_2,k_3\in\mathbb{Z}\\-k_1+k_2+k_3=k}}\fe^{-2i\nu kk_1}(\fe^{2i\nu k_2k_3}-1)k_1\overline{\widehat{w}}_{k_1}(\nu)\widehat{w}_{k_2}(\nu)\widehat{w}_{k_3}(\nu)\,d\nu\,\fe^{ikx}\\
&\quad+i\fe^{i\tau\partial_x^2}\int_0^\tau\sum\limits_{k\in\mathbb{Z}}\sum\limits_{\substack{k_1,k_2,k_3\in\mathbb{Z}\\-k_1+k_2+k_3=k}}\fe^{-2i\nu kk_1}k_1\overline{\widehat{w}}_{k_1}(\nu)\bigl(\widehat{w}_{k_2}(\nu)\widehat{w}_{k_3}(\nu)-\widehat{w}_{k_2}(0)\widehat{w}_{k_3}(0)\bigr)\,d\nu\,\fe^{ikx}\\
&\quad+i\fe^{i\tau\partial_x^2}\int_0^\tau\sum\limits_{k\in\mathbb{Z}}\sum\limits_{\substack{k_1,k_2,k_3\in\mathbb{Z}\\-k_1+k_2+k_3=k}}\fe^{-2i\nu kk_1}k_1\Bigl(\overline{\widehat{w}}_{k_1}(\nu)-\overline{\widehat{w}}_{k_1}(0)-\mathcal{F}_x\bigl(\overline{\phi}(w(0),\nu)\bigr)_{k_1}\Bigr)\\
&\hspace*{8cm}\times\widehat{w}_{k_2}(0)\widehat{w}_{k_3}(0)\,d\nu\,\fe^{ikx}\\
&\quad+i\fe^{i\tau\partial_x^2}\int_0^\tau\sum\limits_{k\in\mathbb{Z}}\sum\limits_{\substack{k_1,k_2,k_3\in\mathbb{Z}\\-k_1+k_2+k_3=k}}\fe^{-2i\nu kk_1}k_1\mathcal{F}_x\bigl(\overline{\phi}(w(0),\nu)\bigr)_{k_1}\widehat{w}_{k_2}(0)\widehat{w}_{k_3}(0)\,d\nu\,\fe^{ikx}\\
&=\mathcal{E}_{1,1}(t_n)+\mathcal{E}_{1,2}(t_n)+\mathcal{E}_{1,3}(t_n)+\mathcal{E}_{1,4}(t_n),
\end{align*}
where 
\begin{equation}\label{defphi}
\begin{aligned}
\phi(w(0),\nu)&=-\frac{i}{2}\partial_x^{-1}g_1(w(0),\nu)\\
&=i\int_0^\nu\sum\limits_{k\in\mathbb{Z}\setminus\{0\}}\sum\limits_{\substack{k_1,k_2,k_3\in\mathbb{Z}\\-k_1+k_2+k_3=k}}\fe^{-2i\xi kk_1}k_1\overline{\widehat{w}}_{k_1}(0)\widehat{w}_{k_2}(0)\widehat{w}_{k_3}(0)\,d\xi\,\fe^{ikx}.
\end{aligned}
\end{equation}

For $\mathcal{E}_{1,1}(t_n)$, we first deduce that
\begin{align*}
\|\mathcal{E}_{1,1}(t_n)\|_s&\lesssim\tau\sup\limits_{\nu\in[0,\tau]}\nu\,\Big\|\sum\limits_{k\in\mathbb{Z}}\sum\limits_{\substack{k_1,k_2,k_3\in\mathbb{Z}\\-k_1+k_2+k_3=k}}k_1k_2k_3\,\varphi_1(2i\nu k_2k_3)\overline{\widehat{w}}_{k_1}(\nu)\widehat{w}_{k_2}(\nu)\widehat{w}_{k_3}(\nu)\fe^{ikx}\Big\|_s.
\end{align*}
Here, we used the fact that, for any functions
$$
f(x)=\sum\limits_{k\in\mathbb{Z}}\widehat{f}_k\fe^{ikx},\qquad g(x)=\sum\limits_{k\in\mathbb{Z}}\widehat{g}_k\fe^{ikx},
$$
the coefficient-wise bound $|\widehat f_k|\le C|\widehat g_k|$ for all $k\in\mathbb Z$ with $C$ independent of $k$ implies
\begin{equation}\label{coeff}
\|f\|_s\lesssim\|g\|_s, \qquad \text{for all}\ s\in\mathbb{R}.
\end{equation}
Noting that $\varphi_1$ is uniformly bounded in $i\mathbb{R}$, by \eqref{uvbound} and \eqref{cbe}, we obtain
\begin{equation}\label{e11}
\begin{aligned}
\|\mathcal{E}_{1,1}(t_n)\|_s&\lesssim\tau\sup\limits_{\nu\in[0,\tau]}\nu\,\Big\|\sum\limits_{k_1,k_2,k_3\in\mathbb{Z}}\bigl(\big|k_1\overline{\widehat{w}}_{k_1}(\nu)\big|\fe^{-ik_1x}\bigr)\bigl(\big|k_2\widehat{w}_{k_2}(\nu)\big|\fe^{ik_2x}\bigr)\bigl(\big|k_3\widehat{w}_{k_3}(\nu)\big|\fe^{ik_3x}\bigr)\Big\|_s\\
&\lesssim\tau^2\sup\limits_{\nu\in[0,\tau]}\Big\|\Bigl(\sum\limits_{k\in\mathbb{Z}}\big|k\widehat{w}_k(\nu)\big|\fe^{ikx}\Bigr)^3\Big\|_s\lesssim\tau^2\sup\limits_{\nu\in[0,\tau]}\|w(\nu)\|^3_{s+1}\lesssim\tau^2.
\end{aligned}
\end{equation}
Similarly, for $\mathcal{E}_{1,2}(t_n)$, we derive that
\begin{align*}
\|\mathcal{E}_{1,2}(t_n)\|_s&\lesssim\tau\sup\limits_{\nu\in[0,\tau]}\Big\|\sum\limits_{k\in\mathbb{Z}}\sum\limits_{\substack{k_1,k_2,k_3\in\mathbb{Z}\\-k_1+k_2+k_3=k}}k_1\overline{\widehat{w}}_{k_1}(\nu)\bigl(\widehat{w}_{k_2}(\nu)\widehat{w}_{k_3}(\nu)-\widehat{w}_{k_2}(0)\widehat{w}_{k_3}(0)\bigr)\fe^{ikx}\Big\|_s\\
&=\tau\sup\limits_{\nu\in[0,\tau]}\Big\|\sum\limits_{k_1,k_2,k_3\in\mathbb{Z}}\bigl(k_1\overline{\widehat{w}}_{k_1}(\nu)\fe^{-ik_1x}\bigr)\bigl(\widehat{w}_{k_2}(\nu)\widehat{w}_{k_3}(\nu)-\widehat{w}_{k_2}(0)\widehat{w}_{k_3}(0)\bigr)\fe^{i(k_2+k_3)x}\Big\|_s\\
&=\tau\sup\limits_{\nu\in[0,\tau]}\Big\|\Bigl(\sum\limits_{k_1\in\mathbb{Z}}k_1\overline{\widehat{w}}_{k_1}(\nu)\fe^{-ik_1x}\Bigr)\bigl(w^2(\nu)-w^2(0)\bigr)\Big\|_s\\
&\lesssim\tau\sup\limits_{\nu\in[0,\tau]}\|w(\nu)\|_{s+1}\|w(\nu)-w(0)\|_s\|w(\nu)+w(0)\|_s.
\end{align*}
Thus, by \eqref{wdiff}, we obtain
\begin{equation}\label{e12}
\|\mathcal{E}_{1,2}(t_n)\|_s\lesssim\tau^2.
\end{equation}
We now estimate $\mathcal{E}_{1,3}(t_n)$. By \eqref{cbe} and Lemma~\ref{lempi0}, we get
\begin{align*}
\|\Pi_0\mathcal{E}_{1,3}(t_n)\|_s=|\Pi_0\mathcal{E}_{1,3}|&\lesssim\tau\sup\limits_{\nu\in[0,\tau]}\big|\Pi_0\bigl(w^2(0)\partial_x(\overline{w}(\nu)-\overline{w}(0)-\overline{\phi}(w(0),\nu)\bigr)\big|\\
&\lesssim\tau\sup\limits_{\nu\in[0,\tau]}\|w(0)\|^2_s\bigl(\|w(\nu)-w(0)\|_{\frac12}+\|\phi(w(0),\nu)\|_\frac12\bigr).
\end{align*}
We have to estimate $\|\phi(w(0),\nu)\|_\frac12$. By \eqref{defphi} and \eqref{coeff}, we have
\begin{equation}\label{phis}
\begin{aligned}
\|\phi(w(0),\nu)\|_\frac12&\lesssim\|\phi(w(0),\nu)\|_s\lesssim\tau\Big\|\sum\limits_{k\in\mathbb{Z}\setminus\{0\}}\sum\limits_{\substack{k_1,k_2,k_3\in\mathbb{Z}\\-k_1+k_2+k_3=k}}k_1\overline{\widehat{w}}_{k_1}(0)\widehat{w}_{k_2}(0)\widehat{w}_{k_3}(0)\fe^{ikx}\Big\|_s\\
&\lesssim\tau\Big\|\sum\limits_{\substack{k_1,k_2,k_3\in\mathbb{Z}\\k_1\neq k_2+k_3}}\bigl(\big|k_1\overline{\widehat{w}}_{k_1}(0)\big|\fe^{-ik_1x}\bigr)\bigl(\big|\widehat{w}_{k_2}(0)\big|\fe^{ik_2x}\bigr)\bigl(\big|\widehat{w}_{k_3}(0)\fe^{ik_3x}\big|\bigr)\Big\|_s\lesssim\tau\|w(0)\|^3_{s+1}.
\end{aligned}
\end{equation}
Thus, by \eqref{wdiff}, we obtain
\begin{equation}\label{e131}
\|\Pi_0\mathcal{E}_{1,3}(t_n)\|_s\lesssim\tau^2.
\end{equation}
Next, by substituting
$$
w(\nu)-w(0)=\int_0^\nu\sum\limits_{j=1}^4\fe^{-i\xi\partial_x^2}F_j\bigl(\fe^{i\xi\partial_x^2}w(\xi)\bigr)d\xi,
$$
we derive that
\begin{align*}
\|(I&-\Pi_0)\mathcal{E}_{1,3}(t_n)\|_s\lesssim\tau^2\sup\limits_{\xi\in[0,\tau]}\big\|(F_2+F_3+F_4)\bigl(\fe^{i\xi\partial_x^2}w(\xi)\bigr)\big\|_{s+1}\\
&+\bigg\|(I-\Pi_0)\int_0^\tau\fe^{-i\nu\partial_x^2}\biggl(\fe^{-i\nu\partial_x^2}\partial_x\Bigl(\int_0^\nu\fe^{i\xi\partial_x^2}\overline{F}_1\bigl(\fe^{i\xi\partial_x^2}w(\xi)\bigr)\,d\xi-\overline{\phi}(w(0),\nu)\Bigr)\fe^{i\nu\partial_x^2}w^2(0)\biggr)\,d\nu\bigg\|_s.
\end{align*}
By \eqref{cbe} we directly get
$$
\sup\limits_{\xi\in[0,\tau]}\big\|(F_2+F_3+F_4)\bigl(\fe^{i\xi\partial_x^2}w(\xi)\bigr)\big\|_{s+1}\lesssim1.
$$
Moreover, calling
$$
\overline{f}(\nu)=(I-\Pi_0)\int_0^\nu\fe^{i\xi\partial_x^2}\overline{F}_1\bigl(\fe^{i\xi\partial_x^2}w(\xi)\bigr)\,d\xi-\overline{\phi}(w(0),\nu),
\qquad g(\nu)=w^2(0),
$$
and noting that $\partial_x=\partial_x(I-\Pi_0)$ and $f(0)=\partial_\nu g=0$, by \eqref{cbe}, \eqref{cbepm} and Lemma~\ref{lemintpart}, we deduce that
\begin{align*}
\bigg\|(I&-\Pi_0)\int_0^\tau\fe^{-i\nu\partial_x^2}\biggl(\fe^{-i\nu\partial_x^2}\partial_x\Bigl(\int_0^\nu\fe^{i\xi\partial_x^2}\overline{F}_1\bigl(\fe^{i\xi\partial_x^2}w(\xi)\bigr)\,d\xi-\overline{\phi}(w(0),\nu)\Bigr)\fe^{i\nu\partial_x^2}w^2(0)\biggr)\,d\nu\bigg\|_s\\
&\lesssim\Big\|(I-\Pi_0)\int_0^\tau\fe^{-i\nu\partial_x^2}F_1\bigl(\fe^{i\nu\partial_x^2}w(\nu)\bigr)\,d\nu-\phi(w(0),\tau)\Big\|_{s-1}\|w(0)\|^2_s\\
&\quad+\tau\sup\limits_{\nu\in[0,\tau]}\big\|(I-\Pi_0)\fe^{-i\nu\partial_x^2}F_1\bigl(\fe^{i\nu\partial_x^2}w(\nu)\bigr)-\partial_\nu\phi(w(0),\nu)\big\|_{s-1}\|w(0)\|^2_s.
\end{align*}
Thus by \eqref{defphi}, we have
\begin{align*}
\Big\|(I&-\Pi_0)\int_0^\tau\fe^{-i\nu\partial_x^2}F_1\bigl(\fe^{i\nu\partial_x^2}w(\nu)\bigr)\,d\nu-\phi(w(0),\tau)\Big\|_{s-1}\\
&\le\Big\|\int_0^\tau\sum\limits_{k\in\mathbb{Z}\setminus\{0\}}\sum\limits_{\substack{k_1,k_2,k_3\in\mathbb{Z}\\-k_1+k_2+k_3=k}}\fe^{-2i\nu kk_1}(\fe^{2i\nu k_2k_3}-1)k_1\overline{\widehat{w}}_{k_1}(\nu)\widehat{w}_{k_2}(\nu)\widehat{w}_{k_3}(\nu)\,d\nu\,\fe^{ikx}\Big\|_{s-1}\\
&\quad+\Big\|\int_0^\tau\!\sum\limits_{k\in\mathbb{Z}\setminus\{0\}}\!\sum\limits_{\substack{k_1,k_2,k_3\in\mathbb{Z}\\-k_1+k_2+k_3=k}}\!\fe^{-2i\nu kk_1}k_1\overline{\widehat{w}}_{k_1}(\nu)\bigl(\widehat{w}_{k_2}(\nu)\widehat{w}_{k_3}(\nu)-\widehat{w}_{k_2}(0)\widehat{w}_{k_3}(0)\bigr)\,d\nu\,\fe^{ikx}\Big\|_{s-1}\\
&\quad+\Big\|\int_0^\tau\!\sum\limits_{k\in\mathbb{Z}\setminus\{0\}}\!\sum\limits_{\substack{k_1,k_2,k_3\in\mathbb{Z}\\-k_1+k_2+k_3=k}}\!\fe^{-2i\nu kk_1}k_1\bigl(\overline{\widehat{w}}_{k_1}(\nu)-\overline{\widehat{w}}_{k_1}(0)\bigr)\widehat{w}_{k_2}(0)\widehat{w}_{k_3}(0)\bigr)\,d\nu\,\fe^{ikx}\Big\|_{s-1}.
\end{align*}
Note that the first and second terms in the above estimate can be handled by a similar argument as given for $\mathcal{E}_{1,1}(t_n)$ and $\mathcal{E}_{1,2}(t_n)$, respectively. Moreover, for the third term, by transforming back to the physical space, by \eqref{cbe}, \eqref{cbepm} and \eqref{coeff}, we obtain
\begin{align*}
\Big\|\int_0^\tau&\sum\limits_{k\in\mathbb{Z}\setminus\{0\}}\sum\limits_{\substack{k_1,k_2,k_3\in\mathbb{Z}\\-k_1+k_2+k_3=k}}\fe^{-2i\nu kk_1}k_1\bigl(\overline{\widehat{w}}_{k_1}(\nu)-\overline{\widehat{w}}_{k_1}(0)\bigr)\widehat{w}_{k_2}(0)\widehat{w}_{k_3}(0)\bigr)\,d\nu\,\fe^{ikx}\Big\|_{s-1}\\
&\lesssim\tau\sup_{\nu\in[0,\tau]}\big\|\partial_x\bigl(\overline{w}(\nu)-\overline{w}(0)\bigr)w^2(0)\big\|_{s-1}\lesssim\tau\|w(\nu)-w(0)\|_s\|w(0)\|^2_s\lesssim\tau^2.
\end{align*}
This shows that
$$
\Big\|(I-\Pi_0)\int_0^\tau\fe^{-i\nu\partial_x^2}F_1\bigl(\fe^{i\nu\partial_x^2}w(\nu)\bigr)\,d\nu-\phi(w(0),\tau)\Big\|_{s-1}\lesssim\tau^2.
$$
Next, by \eqref{defphi}, we derive that
$$
\phi_\nu(w(0),\nu)=i\sum\limits_{k\in\mathbb{Z}\setminus\{0\}}\sum\limits_{\substack{k_1,k_2,k_3\in\mathbb{Z}\\-k_1+k_2+k_3=k}}\fe^{-2i\nu kk_1}k_1\overline{\widehat{w}}_{k_1}(0)\widehat{w}_{k_2}(0)\widehat{w}_{k_3}(0)\,\fe^{ikx}.
$$
Therefore, by a similar argument, we obtain
\begin{equation}\label{e133}
\begin{aligned}
\sup\limits_{\nu\in[0,\tau]}&\big\|(I-\Pi_0)\fe^{-i\nu\partial_x^2}F_1\bigl(\fe^{i\nu\partial_x^2}w(\nu)\bigr)-\partial_\nu\phi\bigl(w(0),\nu\bigr)\big\|_{s-1}\\
&\lesssim\sup\limits_{\nu\in[0,\tau]}\Big\|\sum\limits_{k\in\mathbb{Z}\setminus\{0\}}\sum\limits_{\substack{k_1,k_2,k_3\in\mathbb{Z}\\-k_1+k_2+k_3=k}}\fe^{-2i\nu kk_1}(\fe^{2i\nu k_2k_3}-1)k_1\overline{\widehat{w}}_{k_1}(\nu)\widehat{w}_{k_2}(\nu)\widehat{w}_{k_3}(\nu)\,d\nu\,\fe^{ikx}\Big\|_{s-1}\\
&\quad+\sup\limits_{\nu\in[0,\tau]}\Big\|\sum\limits_{k\in\mathbb{Z}\setminus\{0\}}\!\sum\limits_{\substack{k_1,k_2,k_3\in\mathbb{Z}\\-k_1+k_2+k_3=k}}\!\fe^{-2i\nu kk_1}k_1\overline{\widehat{w}}_{k_1}(\nu)\bigl(\widehat{w}_{k_2}(\nu)\widehat{w}_{k_3}(\nu)-\widehat{w}_{k_2}(0)\widehat{w}_{k_3}(0)\bigr)\,\fe^{ikx}\Big\|_{s-1}\\
&\quad+\sup\limits_{\nu\in[0,\tau]}\Big\|\sum\limits_{k\in\mathbb{Z}\setminus\{0\}}\!\sum\limits_{\substack{k_1,k_2,k_3\in\mathbb{Z}\\-k_1+k_2+k_3=k}}\!\fe^{-2i\nu kk_1}k_1\bigl(\overline{\widehat{w}}_{k_1}(\nu)-\overline{\widehat{w}}_{k_1}(0)\bigr)\widehat{w}_{k_2}(0)\widehat{w}_{k_3}(0)\bigr)\,\fe^{ikx}\Big\|_{s-1}\\
&\lesssim\tau.
\end{aligned}
\end{equation}
Thus by combing the above estimates, we arrive at
\begin{equation}\label{e13}
\|\mathcal{E}_{1,3}(t_n)\|_s\lesssim\tau^2.
\end{equation}

Note that $\mathcal{E}_{1,4}(t_n)$ can not be bounded by $\tau^2$ at low regularity. Nevertheless, in the following, we will prove that
\begin{equation}\label{e142}
\|\mathcal{E}_{1,4}(t_n)+\mathcal{E}_2(t_n)\|_s\lesssim\tau^2.
\end{equation}
This also explains why the method \eqref{meth1} approximates the quintic term $F_2$ in a different way as in~\cite{ostfocm}.

We first calculate $\mathcal{E}_{1,4}(t_n)$. By \eqref{defphi}, we have
\begin{equation}\label{e14cal}
\begin{aligned}
\mathcal{E}_{1,4}(t_n)&=i\fe^{i\tau\partial_x^2}\int_0^\tau\sum\limits_{k\in\mathbb{Z}}\sum\limits_{\substack{k_1,k_2,k_3\in\mathbb{Z}\\-k_1+k_2+k_3=k\\k_1\neq0}}\fe^{-2i\nu kk_1}k_1\widehat{w}_{k_2}(0)\widehat{w}_{k_3}(0)\\
&\qquad\qquad\times\frac{i}2\Bigl(\frac{1}{-ik_1}\Bigr)\sum\limits_{\substack{k_4,k_5,k_6\in\mathbb{Z}\\-k_4+k_5+k_6=k_1}}(\fe^{2i\nu k_1k_4}-1)\widehat{w}_{k_4}(0)\overline{\widehat{w}}_{k_5}(0)\overline{\widehat{w}}_{k_6}(0)\,d\nu\,\fe^{ikx}\\
&=-\frac{i}2\fe^{i\tau\partial_x^2}\int_0^\tau\sum\limits_{k\in\mathbb{Z}}\sum\limits_\mathcal{K}(\fe^{-2i\nu(-k_4+k_5+k_6)(k-k_4)}-\fe^{-2i\nu k(-k_4+k_5+k_6)})\\
&\qquad\qquad\times\widehat{w}_{k_2}(0)\widehat{w}_{k_3}(0)\widehat{w}_{k_4}(0)\overline{\widehat{w}}_{k_5}(0)\overline{\widehat{w}}_{k_6}(0)\,d\nu\,\fe^{ikx},
\end{aligned}
\end{equation} 
where
$$
\mathcal{K}=\{(k_2,k_3,k_4,k_5,k_6)\mid k_2+k_3+k_4-k_5-k_6=k\}.
$$
We then calculate $\mathcal{E}_2(t_n)$. Specifically, we have
\begin{align*}
\mathcal{E}_2(t_n)&=\fe^{i\tau\partial_x^2}\int_0^\tau\fe^{-i\nu\partial_x^2}\Bigl(F_2\bigl(\fe^{i\nu\partial_x^2}w(\nu)\bigr)-F_2\bigl(\fe^{i\nu\partial_x^2}w(0)\bigr)\Bigr)\,d\nu\\
&\quad+\frac{i}2\fe^{i\tau\partial_x^2}\int_0^\tau\sum\limits_{k\in\mathbb{Z}}\sum\limits_\mathcal{K}\bigl(\fe^{i\nu(k^2+k_5^2+k_6^2-k_2^2-k_3^2-k_4^2)}-\fe^{2i\nu(k_4-k_5-k_6)^2}\bigr)\\
&\qquad\qquad\qquad\qquad\times\widehat{w}_{k_2}(0)\widehat{w}_{k_3}(0)\widehat{w}_{k_4}(0)\overline{\widehat{w}}_{k_5}(0)\overline{\widehat{w}}_{k_6}(0)\,d\nu\,\fe^{ikx}\\
&=\mathcal{E}_{2,1}(t_n)+\mathcal{E}_{2,2}(t_n).
\end{align*}
For $\mathcal{E}_{2,1}(t_n)$, by \eqref{cbe} and \eqref{wdiff}, we obtain
\begin{equation}\label{e21}
\|\mathcal{E}_{2,1}(t_n)\|_s\lesssim\tau\sup\limits_{\nu\in[0,\tau]}\|w(\nu)-w(0)\|_s\sum\limits_{j=0}^4\|w(\nu)\|^j_s\|w(0)\|^{4-j}_s\lesssim\tau^2.
\end{equation}

We next combine $\mathcal{E}_{1,4}(t_n)$ and $\mathcal{E}_{2,2}(t_n)$. Indeed, we have
\begin{align*}
\mathcal{E}_{1,4}(t_n)&+\mathcal{E}_{2,2}(t_n)=\frac{i}2\fe^{i\tau\partial_x^2}\int_0^\tau\sum\limits_{k\in\mathbb{Z}}\sum\limits_\mathcal{K}\Bigl(\bigl(\fe^{i\nu\rho_1(k)}-\fe^{i\nu\rho_2(k)}\bigr)+\bigl(\fe^{i\nu\rho_3(k)}-\fe^{i\nu\rho_4(k)}\bigr)\Bigr)\\
&\qquad\qquad\qquad\qquad\times\widehat{w}_{k_2}(0)\widehat{w}_{k_3}(0)\widehat{w}_{k_4}(0)\overline{\widehat{w}}_{k_5}(0)\overline{\widehat{w}}_{k_6}(0)\,d\nu\,\fe^{ikx},
\end{align*}
where 
\begin{alignat*}{2}
&\rho_1(k)=k^2+k_5^2+k_6^2-k_2^2-k_3^2-k_4^2,\qquad &&\rho_2(k)=2(k_4-k_5-k_6)(k_2+k_3-k_5-k_6),\\
&\rho_3(k)=2k(k_4-k_5-k_6), &&\rho_4(k)=2(k_4-k_5-k_6)^2.
\end{alignat*}
Note that we have
$$
|\rho_1(k)-\rho_2(k)|+|\rho_3(k)-\rho_4(k)|\lesssim\prod\limits_{i=2}^6(1+|k_i|).
$$
Thus, similar to \cite{ostfocm}, we obtain
\begin{equation}\label{e22}
\|\mathcal{E}_{1,4}(t_n)+\mathcal{E}_{2,2}(t_n)\|_s\lesssim\tau\sup_{\nu\in[0,\tau]}\nu\Big\|\prod\limits_{j=2}^6(1+|k_j|)|\widehat{w}_{k_j}(0)|\fe^{ik_jx}\Big\|_s\lesssim\tau^2\|w(0)\|^5_{s+1}. 
\end{equation}
This finally shows \eqref{e142}.

Since $\mu$ can be bounded by \eqref{uvbound}, we shall omit the estimate for $\mathcal{E}_3(t_n)$. The argument is very similar to that given in \cite{ostfocm}. Indeed, by \eqref{wdiff}, we derive that
\begin{equation}\label{e3}
\|\mathcal{E}_3(t_n)\|_s\lesssim\tau^2.
\end{equation}

We finally estimate $\mathcal{E}_4(t_n)$. Since $\psi$ is independent of $x$, we actually have
$$
\fe^{-i\nu\partial_x^2}\bigl(\psi(v(t_n+\nu))v(t_n+\nu)\bigr)=\psi(v(t_n+\nu))\fe^{-i\nu\partial_x^2}v(t_n+\nu).
$$
Thus by \eqref{wnuandv} and \eqref{wdiff}, we deduce that
\begin{align*}
\|\mathcal{E}_4(t_n)\|_s&\lesssim\Big\|\int_0^\tau|\psi(v(t_n+\nu))-\psi(v(t_n))|w(\nu)\,d\nu\Big\|_s+\tau\sup_{\nu\in[0,\tau]}|\psi(v(t_n))|\|w(\nu)-w(0)\|_s\\
&\lesssim\tau\sup_{\nu\in[0,\tau]}\big|\psi(w(\nu))-\psi(w(0))\big|\|w(\nu)\|_s+\tau^2|\psi(w(0))|
\end{align*}
since $\psi(v(t_n+\nu))=\psi(w(\nu))$. Moreover, combining \eqref{uvbound} and \eqref{wdiff} with the Sobolev embedding $H^{\frac14}\hookrightarrow L^4$ and Lemma~\ref{lempi0}, we obtain
\begin{align*}
\big|\psi(w(\nu))-\psi(w(0))\big|\lesssim\|w(\nu)&-w(0)\|_\frac12\bigl(\|w(\nu)\|_\frac12+\|w(0)\|_\frac12\bigr)\\
&+\|w(\nu)-w(0)\|_\frac14\bigl(\|w(\nu)\|_\frac14+\|w(0)\|_\frac14\bigr)^3\lesssim\tau,
\end{align*}
and
$$
\big|\psi(w(0))\big|\lesssim\|w(0)\|^2_\frac12+\|w(0)\|^4_\frac14\lesssim1.
$$
This yields
\begin{equation}\label{e4}
\|\mathcal{E}_4(t_n)\|_s\lesssim\tau^2.
\end{equation}

We conclude the proof of Proposition~\ref{proplocal1} by collecting \eqref{e11}, \eqref{e12}, \eqref{e13}, \eqref{e142}, \eqref{e3} and~\eqref{e4}.
\end{proof}

\section{Global error analysis}\label{sectionglobal}

In this section, we will estimate the global error of method \eqref{meth1} and prove the estimate \eqref{estorder1} given in Theorem~\ref{mainthm}.
\begin{proof}[Proof of Theorem~\ref{mainthm}]
Let $e_n=v_n-v(t_n)$ denote the global error at time $t_n$. With the help of Proposition~\ref{proplocal1}, we get
\begin{equation}\label{glo1st}
\begin{aligned}
\|e_n\|_s=\|\Psi_1^\tau(v_{n-1})&-v(t_n)\|_s\leq\|\Psi_1^\tau(v_{n-1})-\Psi_1^\tau(v(t_{n-1}))\|_s+\|\Psi_1^\tau(v(t_{n-1}))-v(t_n)\|_s\\
&\le\|\Psi_1^\tau(v(t_{n-1})+e_{n-1})-\Psi_1^\tau(v(t_{n-1}))\|_s+C\tau^2.
\end{aligned}
\end{equation}
From \eqref{meth1} and the fact that $g_0$ is a linear operator in its first variable, it follows that
\begin{equation}\label{glo1all}
\begin{aligned}
\|\Psi_1^\tau&(v(t_{n-1})+e_{n-1})-\Psi_1^\tau(v(t_{n-1}))\|_s\le\big\|e_{n-1}-\frac{i}2\partial_x^{-1}g_0\bigl(e_{n-1},v(t_{n-1}),\tau\bigr)\big\|_s\\
&+\frac12\big\|g_0\bigl(v(t_{n-1})+e_{n-1},v(t_{n-1})+e_{n-1},\tau\bigr)-g_0\bigl(v(t_{n-1})+e_{n-1},v(t_{n-1}),\tau\bigr)\|_{s-1}\\
&+\tau\big|\Pi_0\bigl((v(t_{n-1})+e_{n-1})^2\partial_x(\overline{v}(t_{n-1})+\overline{e}_{n-1})-v^2(t_{n-1})\overline{v}_x(t_{n-1})\bigr)\big|\\
&+\frac12\tau\big\|(v(t_{n-1})+e_{n-1})^2\varphi_1(-2i\tau\partial_x^2)\bigl((\overline{v}(t_{n-1})+\overline{e}_{n-1})^2(v(t_{n-1})+e_{n-1})\bigr)\\
&\qquad\qquad\qquad-v^2(t_{n-1})\varphi_1(-2i\tau\partial_x^2)\bigl(\overline{v}^2(t_{n-1})v(t_{n-1})\bigr)\big\|_s\\
&+\mu\tau\|(v(t_{n-1})+e_{n-1})^2\varphi_1(-2i\tau\partial_x^2)(\overline{v}(t_{n-1})+\overline{e}_{n-1})-v^2(t_{n-1})\varphi_1(-2i\tau\partial_x^2)\overline{v}(t_{n-1})\|_s\\
&+\tau\big|\psi(v(t_{n-1})+e_{n-1})-\psi(v(t_{n-1}))\big|\bigl(\|v(t_{n-1})\|_s+\|e_{n-1}\|_s\bigr)+\tau|\psi(v(t_{n-1}))|\|e_{n-1}\|_s.
\end{aligned}
\end{equation}
By \eqref{uvbound}, \eqref{cbe} and Lemma~\ref{lempi0}, we derive that
\begin{align*}
\big|\Pi_0&\bigl((v(t_{n-1})+e_{n-1})^2\partial_x(\overline{v}(t_{n-1})+\overline{e}_{n-1})-v^2(t_{n-1})\partial_x\overline{v}(t_{n-1})\bigr)\big|\\
&\lesssim\|v(t_{n-1})\|^2_s\|e_{n-1}\|_s+\|v(t_{n-1})\|_s\|e_{n-1}\|^2_s+\|e_{n-1}\|^3_s\lesssim\|e_{n-1}\|_s+\|e_{n-1}\|^2_s+\|e_{n-1}\|^3_s.
\end{align*}
Similarly, we have
\begin{align*}
\big|\psi&(v(t_{n-1})+e_{n-1})-\psi(v(t_{n-1}))\big|\\
&\lesssim\|v(t_{n-1})\|_s\|e_{n-1}\|_s+\|e_{n-1}\|^2_s+\sum\limits_{j=1}^4\|v(t_{n-1})\|^{4-j}_s\|e_{n-1}\|^j_s\lesssim\sum\limits_{j=1}^4\|e_{n-1}\|^j_s
\end{align*}
and
$$
|\psi(v(t_{n-1}))|\lesssim\|v(t_{n-1})\|_s^2+\|v(t_{n-1})\|_s^4\lesssim1.
$$
Moreover, for the quintic and cubic terms, we get a similar estimate as that given in \cite{ostfocm}
\begin{align*}
\big\|(v(t_{n-1})+e_{n-1})^2&\varphi_1(-2i\tau\partial_x^2)\bigl((\overline{v}(t_{n-1})+\overline{e}_{n-1})^2(v(t_{n-1})+e_{n-1})\bigr)\\
&-v^2(t_{n-1})\varphi_1(-2i\tau\partial_x^2)\bigl(\overline{v}^2(t_{n-1})v(t_{n-1})\bigr)\big\|_s\lesssim\sum\limits_{j=1}^5\|e_{n-1}\|^j_s
\end{align*}
and
$$
\|(v(t_{n-1})+e_{n-1})^2\varphi_1(-2i\tau\partial_x^2)(\overline{v}(t_{n-1})+\overline{e}_{n-1})-v^2(t_{n-1})\varphi_1(-2i\tau\partial_x^2)\overline{v}(t_{n-1})\|_s\lesssim\sum\limits_{j=1}^3\|e_{n-1}\|^j_s.
$$
Therefore, we already have the estimates for all the terms in \eqref{glo1all} except for those containing $g_0$, i.e.,
\begin{equation}\label{glo1rest}
\begin{aligned}
\|e_n\|_s&\le\frac12\big\|g_0\bigl(v(t_{n-1})+e_{n-1},v(t_{n-1})+e_{n-1},\tau\bigr)-g_0\bigl(v(t_{n-1})+e_{n-1},v(t_{n-1}),\tau\bigr)\big\|_{s-1}\\
&\quad+\big\|e_{n-1}-\frac{i}2\partial_x^{-1}g_0\bigl(e_{n-1},v(t_{n-1}),\tau\bigr)\big\|_s+C\tau\sum\limits_{j=1}^5\|e_{n-1}\|^j_s+C\tau^2.
\end{aligned}
\end{equation}

We now estimate the terms containing $g_0$. Since $g_0$ is a linear operator in the first variable, we can estimate $g_0\bigl(v(t_{n-1}),\cdot,\tau\bigr)$ and $g_0\bigl(e_{n-1},\cdot,\tau\bigr)$ separately to avoid the appearance of $\|e_{n-1}\|_{s+1}$ in the estimate. By \eqref{g0int}, \eqref{cbe} and \eqref{coeff}, we obtain
\begin{equation}\label{g0vpart}
\begin{aligned}
\|g_0\bigl(v(t_{n-1}),&v(t_{n-1})+e_{n-1},\tau\bigr)-g_0\bigl(v(t_{n-1}),v(t_{n-1}),\tau\bigr)\big\|_{s-1}\\
&\lesssim\tau\big\|\overline{v}_x(t_{n-1})\bigl((v(t_{n-1})+e_{n-1})^2-v^2(t_{n-1})\bigr)\big\|_s\lesssim\tau\|e_{n-1}\|_s+\tau\|e_{n-1}\|^2_s.
\end{aligned}
\end{equation}
Moreover, by Young's inequality, \eqref{g0diff}, \eqref{g0int}, \eqref{cbe} and \eqref{cbepm}, we deduce that
\begin{equation}\label{g0epart}
\begin{aligned}
\|g_0&\bigl(e_{n-1},v(t_{n-1})+e_{n-1},\tau\bigr)-g_0\bigl(e_{n-1},v(t_{n-1}),\tau\bigr)\big\|_{s-1}\\
&\lesssim\|g_0\bigl(e_{n-1},v(t_{n-1})+e_{n-1},\tau\bigr)-g_0\bigl(e_{n-1},v(t_{n-1}),\tau\bigr)\big\|^\frac12_{s-2}\\
&\quad\times\|g_0\bigl(e_{n-1},v(t_{n-1})+e_{n-1},\tau\bigr)-g_0\bigl(e_{n-1},v(t_{n-1}),\tau\bigr)\big\|^\frac12_s\\
&\lesssim\tau^\frac12\big\|\partial_xe_{n-1}\bigl((v(t_{n-1})+e_{n-1})^2-v^2(t_{n-1})\bigr)\big\|^\frac12_{s-1}\big\|e_{n-1}\bigl((v(t_{n-1})+e_{n-1})^2-v^2(t_{n-1})\bigr)\big\|^\frac12_s\\
&\lesssim\tau^\frac12\bigl(\|e_{n-1}\|_s^2\|v(t_{n-1})\|_s+\|e_{n-1}\|_s^3\bigr)\lesssim\tau^\frac12\|e_{n-1}\|_s^2+\tau^\frac12\|e_{n-1}\|_s^3.
\end{aligned}
\end{equation}

It remains to estimate $g_0\bigl(e_{n-1},v(t_{n-1}),\tau\bigr)$. In the following, we will prove that
\begin{equation}\label{epn}
\begin{aligned}
\big\|e_{n-1}&-\frac{i}2\partial_x^{-1}g_0\bigl(e_{n-1},v(t_{n-1}),\tau\bigr)\big\|^2_s\\
&\le\|e_{n-1}\|_s^2+\frac14\big\|g_0\bigl(e_{n-1},v(t_{n-1}),\tau\bigr)\big\|^2_{s-1}+\big|\big\langle\langle\partial_x\rangle^se_{n-1},\langle\partial_x\rangle^s\partial_x^{-1}g_0\bigl(e_{n-1},v(t_{n-1}),\tau\bigr)\big\rangle\big|\\
&\le(1+C\tau)\|e_{n-1}\|^2_s,
\end{aligned}
\end{equation}
where $\langle\partial_x\rangle=(I-\partial_x^2)^\frac12.$
First of all, similar to \eqref{g0epart}, we have
\begin{equation}\label{estp}
\begin{aligned}
\big\|g_0&\bigl(e_{n-1},v(t_{n-1}),\tau\bigr)\big\|^2_{s-1}\lesssim\big\|g_0\bigl(e_{n-1},v(t_{n-1}),\tau\bigr)\big\|_{s-2}\big\|g_0\bigl(e_{n-1},v(t_{n-1}),\tau\bigr)\big\|_s\\
&\lesssim\tau\|v^2(t_{n-1})\partial_x\overline{e}_{n-1}\|_{s-1}\|v^2(t_{n-1})\overline{e}_{n-1}\|_s\lesssim\tau\|e_{n-1}\|^2_s.
\end{aligned}
\end{equation}
Next, we will employ the estimate
\begin{equation}\label{wuest}
|\langle\langle\partial_x\rangle^sf,\langle\partial_x\rangle^s\partial_x(\overline{f}g)\rangle|\lesssim\|f\|_s^2\|g\|_{s+1}, \qquad s>\frac12,
\end{equation}
which is proved in exactly the same way as \cite[Lemma~3.3(i)]{wuintpart}. Using \eqref{g0int}, \eqref{wuest} and finally \eqref{cbe}, we deduce that
\begin{equation}\label{estep}
\begin{aligned}
\big|\big\langle\langle\partial_x\rangle^s&e_{n-1},\langle\partial_x\rangle^s\partial_x^{-1}g_0\bigl(e_{n-1},v(t_{n-1}),\tau\bigr)\big\rangle\big|\\
&\lesssim\tau\sup_{\nu\in[0,\tau]}\big|\big\langle\langle\partial_x\rangle^se_{n-1},\langle\partial_x\rangle^s\fe^{-i\nu\partial_x^2}\bigl((\fe^{i\nu\partial_x^2}v^2(t_{n-1}))\partial_x\fe^{-i\nu\partial_x^2}\overline{e}_{n-1}\bigr)\big\rangle\big|\\
&\le\tau\sup_{\nu\in[0,\tau]}\big|\big\langle\langle\partial_x\rangle^s\fe^{i\nu\partial_x^2}e_{n-1},\langle\partial_x\rangle^s\partial_x\bigl((\fe^{-i\nu\partial_x^2}\overline{e}_{n-1})\fe^{i\nu\partial_x^2}v^2(t_{n-1})\bigr)\big\rangle\big|\\
&\quad+\tau\sup_{\nu\in[0,\tau]}|\langle\langle\partial_x\rangle^s\fe^{i\nu\partial_x^2}e_{n-1},\langle\partial_x\rangle^s\bigl((\fe^{-i\nu\partial_x^2}\overline{e}_{n-1})\partial_x\fe^{i\nu\partial_x^2}v^2(t_{n-1})\bigr)\big\rangle\big|\\
&\lesssim\tau\|e_{n-1}\|_s^2\|v^2(t_{n-1})\|_{s+1}+\tau\|e_{n-1}\|_s\|e_{n-1}\|_s\|\partial_xv^2(t_{n-1})\|_s\lesssim\tau\|e_{n-1}\|^2_s.
\end{aligned}
\end{equation}

Thus by \eqref{estp} and \eqref{estep} we conclude \eqref{epn}. Combining \eqref{glo1rest}, \eqref{g0vpart}, \eqref{g0epart}, and \eqref{epn}, we finally obtain
$$
\|e_n\|_s\leq\|e_{n-1}\|_s+C\tau\sum_{j=1}^5\|e_{n-1}\|^j_s+C\tau^\frac12\|e_{n-1}\|^2_s+C\tau^\frac12\|e_{n-1}\|_s^3+C\tau^2.
$$
Using a discrete Gronwall lemma together with the a priori assumption
$$
\|e_n\|_s\le C\tau^\frac12
$$
finally yields \eqref{estorder1}. This ends the proof of Theorem~\ref{mainthm}.
\end{proof}

\section{Numerical Experiments}\label{sectionnumexp}

\begin{figure}
\begin{center}
\subfigure{\includegraphics[height=0.35\textwidth,width=0.45\textwidth]{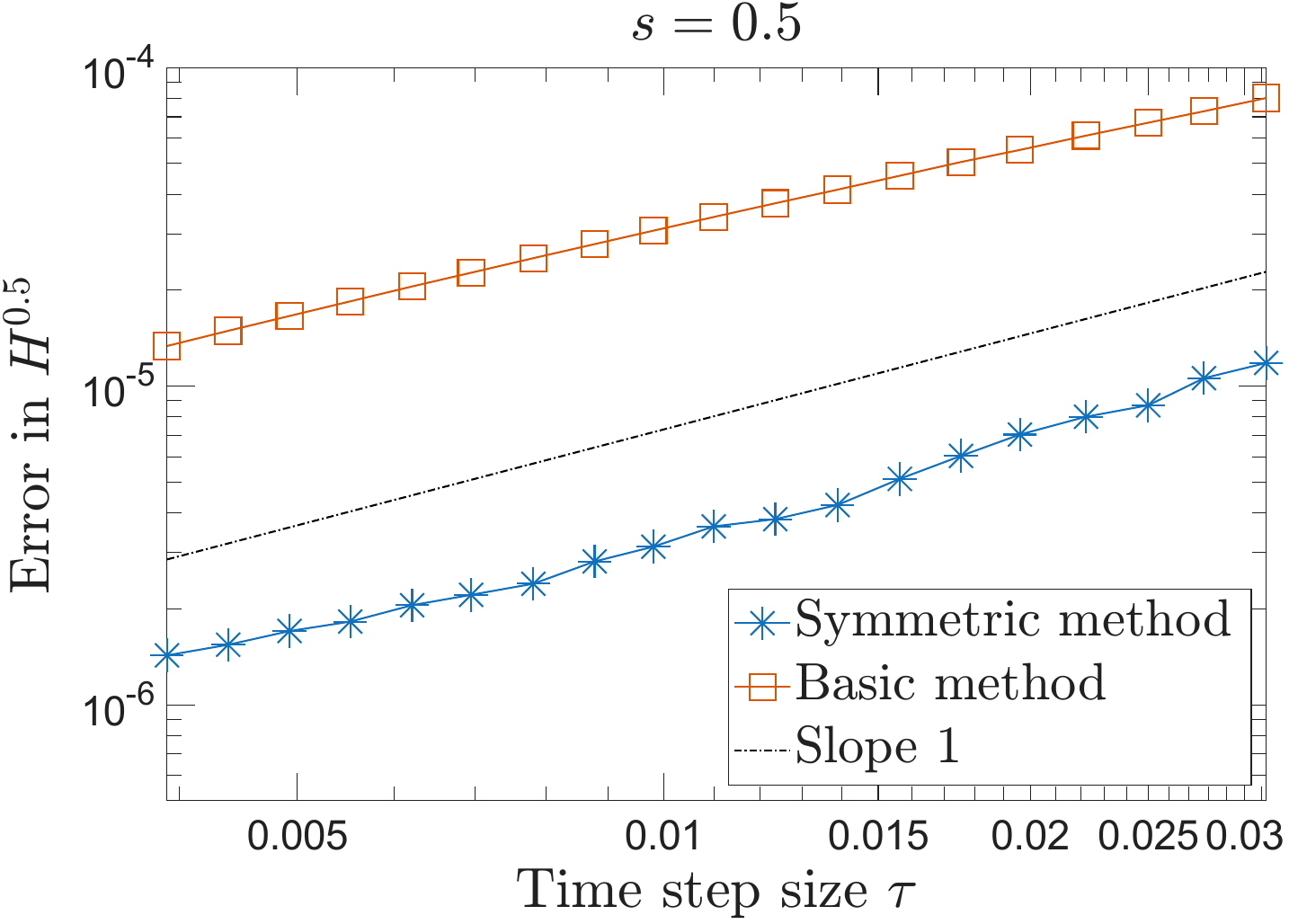}}\hspace{0.4cm}
\subfigure{\includegraphics[height=0.35\textwidth,width=0.45\textwidth]{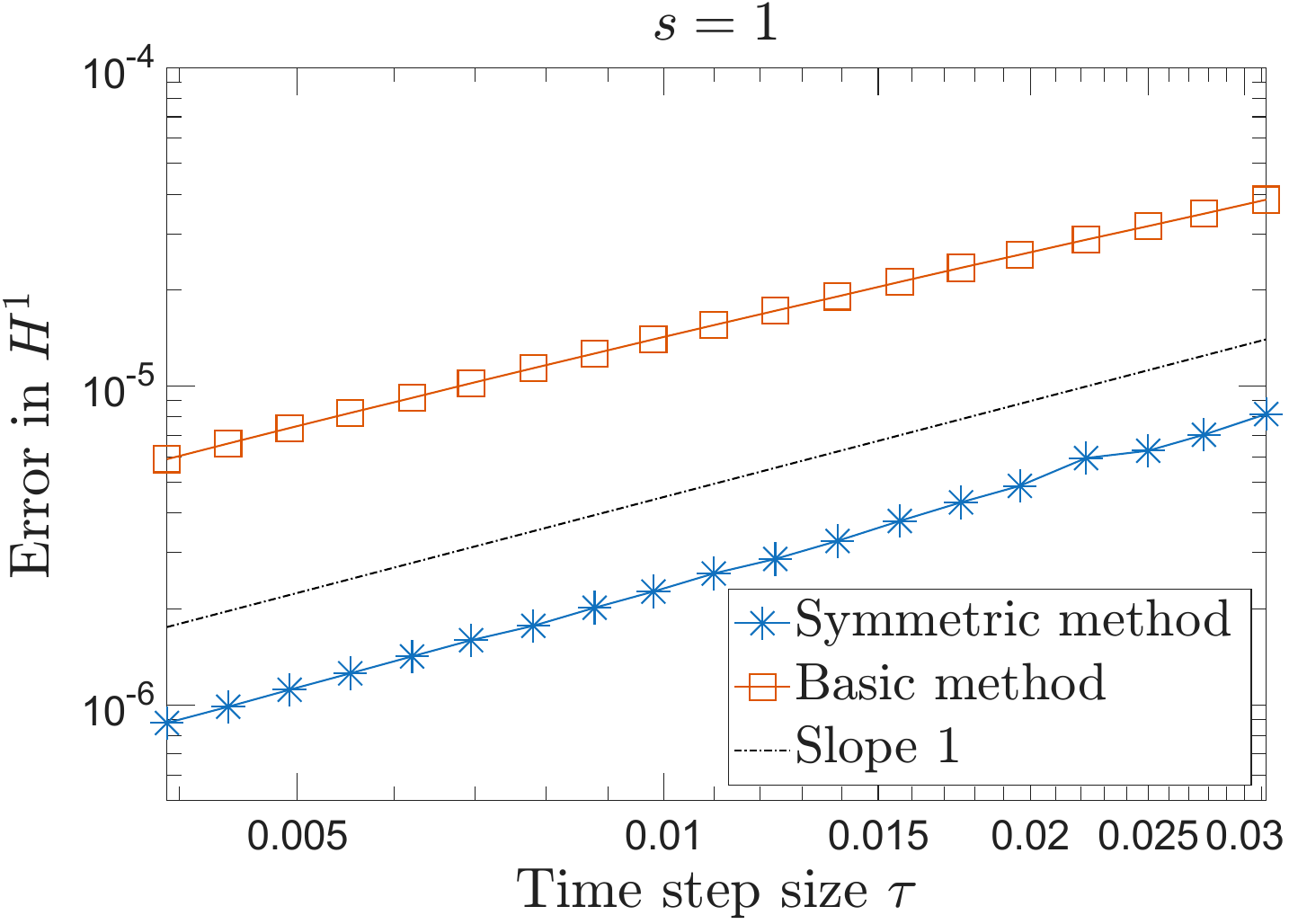}}
\end{center}
\caption{$H^s$ error of our low regularity exponential integrators for rough initial data $u_0\in H^{s+1}$. Left: $s=0.5$; right: $s=1$. \label{figglo}}
\end{figure}

In this section, we provide a numerical illustration of our main result, Theorem ~\ref{mainthm}. We demonstrate the convergence order of our low regularity integrators with rough initial data. Specifically, we take the initial data
$$
u_0(x) = \sum\limits_{k\in\mathbb{Z}}\langle k\rangle^{-(s+1.5+\varepsilon)}\widehat{g}_k \fe^{ikx} \in H^{s+1},
$$
where $\varepsilon$ is a positive constant and the random variables $\widehat{g}_k$ are uniformly distributed in the square~$[-1,1] + i [-1,1]$. In our experiments, we can take $\varepsilon=0$ and choose $s = 0.5,1$. Note that $s=0.5$ is the borderline of our analysis. For spatial discretisation, we use a standard Fourier pseudospectral method (FFT)  and choose the largest Fourier mode as $K=2^{10}$, i.e., the spatial mesh size is set to be equal to $0.0061$. We normalize the $H^{s+1}$ norm of the initial data to $0.5$. As a reference solution, we use our low regularity integrator with $K$ spatial points and a small time step size, $\tau = 2^{-12}$. We set $T=1$ as the final time.

Figure~\ref{figglo} clearly confirms that both methods~\eqref{meth1} and~\eqref{meth2} converge with first order in $H^s$ for initial data $u_0\in H^{s+1}$ (see also Theorem~\ref{mainthm}) for $s=0.5$ and~$1$. Moreover, the global error of both methods is smaller when $s = 1$ than when $s = 0.5$. Additionally, the symmetric method~\eqref{meth2} has a smaller global error than the basic method~\eqref{meth1}.

\begin{figure}
\begin{center}
\subfigure{\includegraphics[height=0.35\textwidth,width=0.45\textwidth]{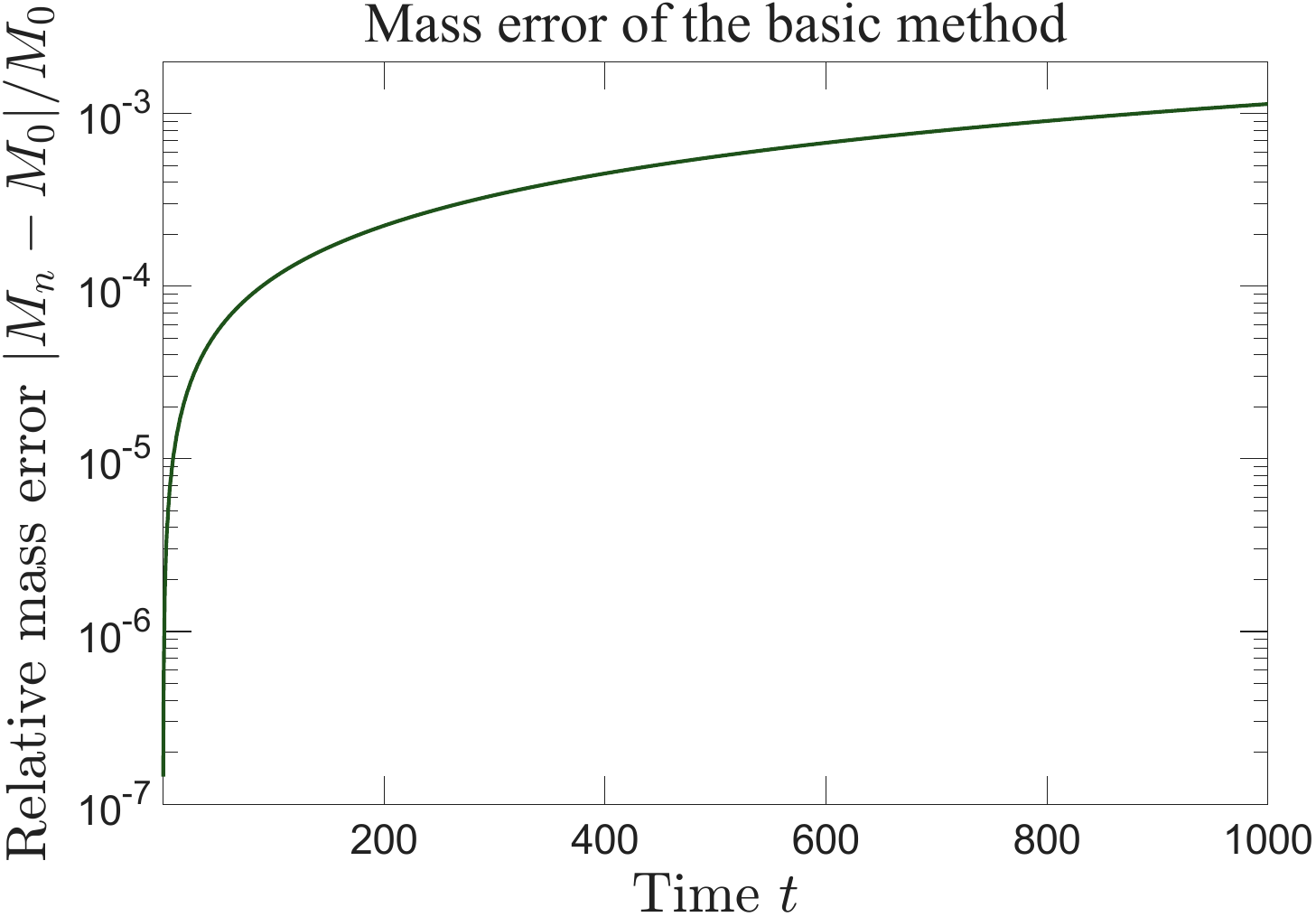}}\hspace{0.4cm}
\subfigure{\includegraphics[height=0.35\textwidth,width=0.45\textwidth]{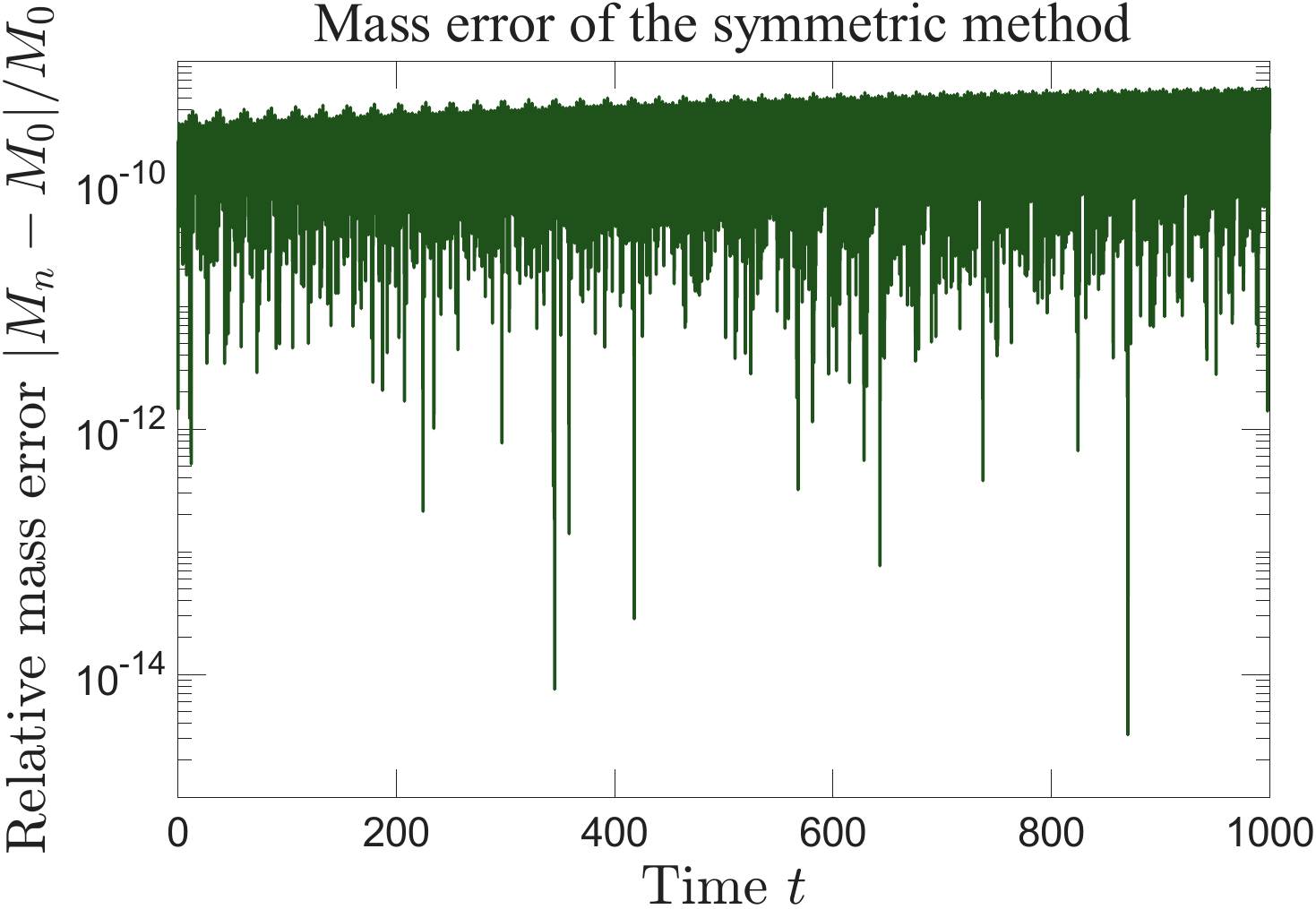}}
\end{center}
\caption{Mass conservation behavior of low regularity exponential integrators for rough initial data $u_0\in H^2$. Left: relative mass error of the basic method \eqref{meth1}; right: relative mass error of the symmetric method \eqref{meth2}.
\label{figmass}}
\end{figure}

Due to the use of the gauge transformation~\eqref{eqgaugedef}, conservation properties, especially mass conservation, are also important for dNLS. Therefore, we also conducted numerical experiments to check the mass~\eqref{mass} and energy~\eqref{hami} conservation of our methods~\eqref{meth1} and~\eqref{meth2}. Through our experiments, we found that $s$, $\tau$, $K$ and $\|u_0\|_{s+1}$ all significantly impact the relative mass and energy error. Here, as an example, we take initial data in $H^2$, choose $K=2^{10}$, $\tau=2^{-10}$, $\|u_0\|_2=0.5$, and the final time to $T=1000$.

Figure~\ref{figmass} shows that both of our methods~\eqref{meth1} and~\eqref{meth2} exhibit good long-term mass behaviour, even at low regularity. In particular, the symmetric method~\eqref{meth2} demonstrates remarkably good mass conservation over long time intervals.

\begin{figure}
\begin{center}
\subfigure{\includegraphics[height=0.35\textwidth,width=0.45\textwidth]{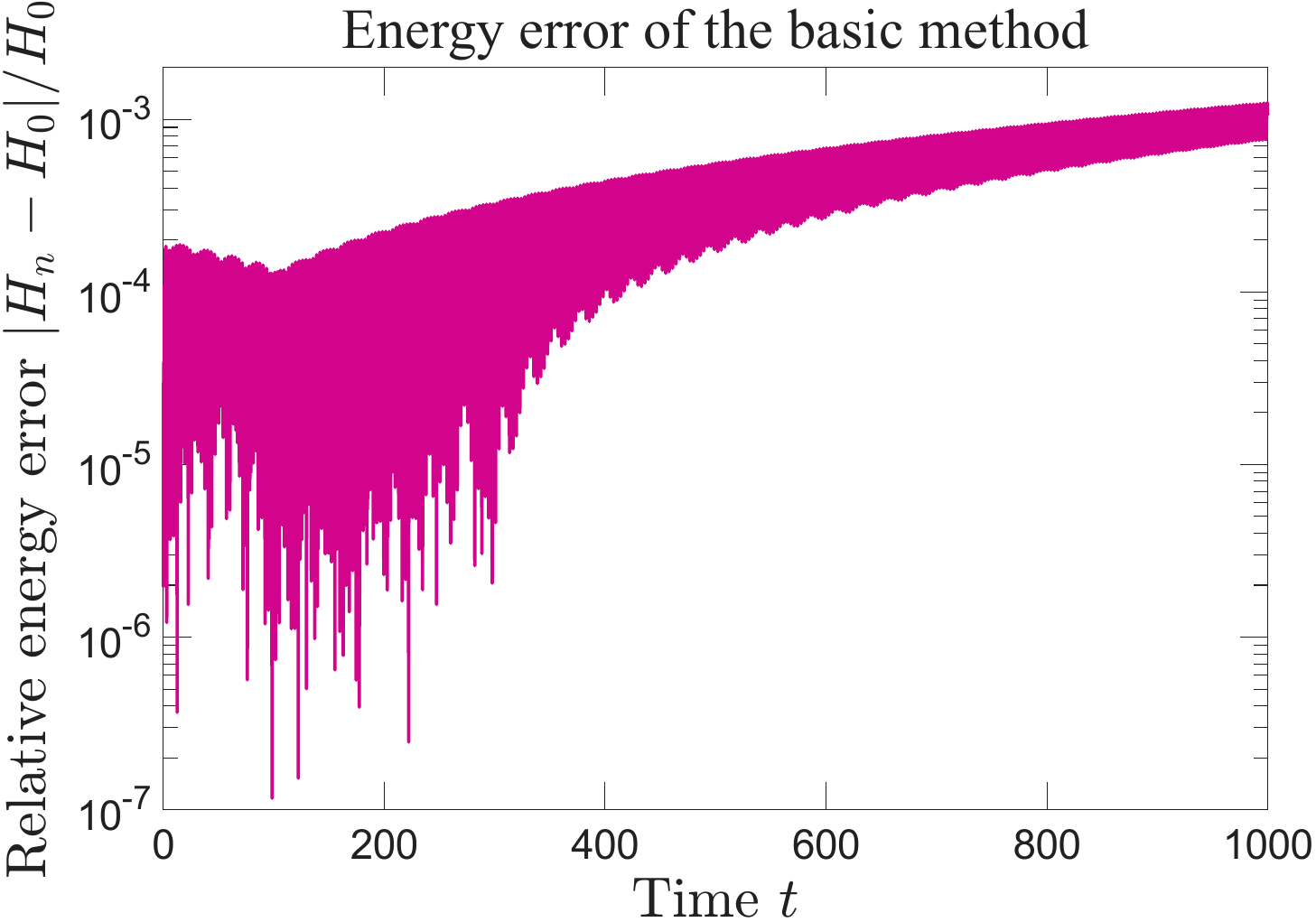}}\hspace{0.4cm}
\subfigure{\includegraphics[height=0.35\textwidth,width=0.45\textwidth]{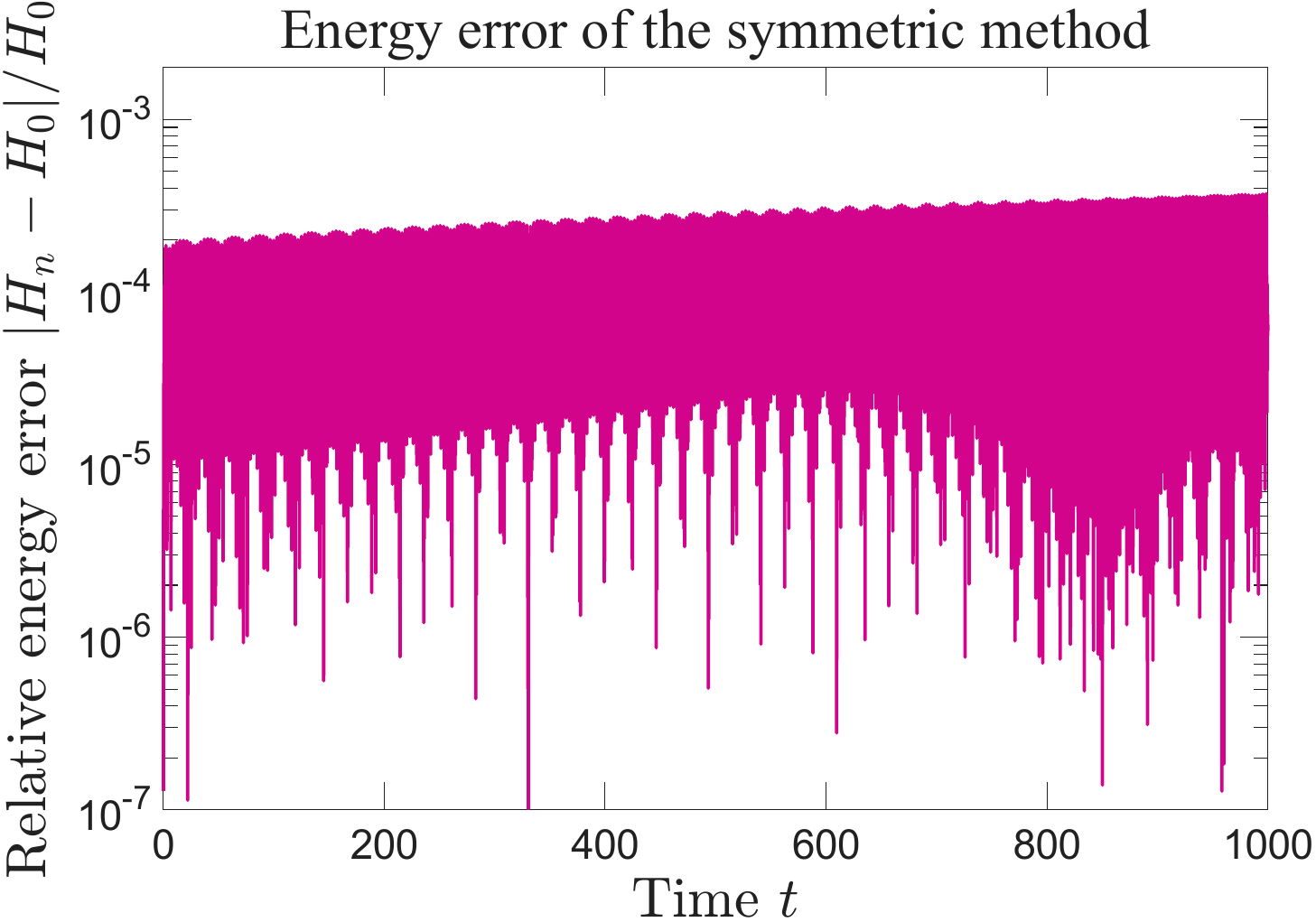}}
\end{center}
\caption{Energy conservation behavior of low regularity exponential integrators for rough initial data $u_0\in H^2$. Left: relative energy error of the basic method \eqref{meth1}; right: relative energy error of the symmetric method \eqref{meth2}.\label{figham}}
\end{figure}

Figure~\ref{figham} illustrates that both of our methods~\eqref{meth1} and~\eqref{meth2} also display good long-term energy behavior at low regularity. Moreover, the symmetric method provides superior energy conservation compared with the basic method.

\end{document}